\documentclass[10pt]{amsart}

\usepackage{amsfonts,amssymb}
\input{diagrams}

\newtheorem{proposition}{Proposition}[section]
\newtheorem{lemma}[proposition]{Lemma}
\newtheorem{corollary}[proposition]{Corollary}
\newtheorem{theorem}[proposition]{Theorem}
\theoremstyle{definition}
\newtheorem{definition}[proposition]{Definition}

\theoremstyle{remark}
\newtheorem{remark}[proposition]{Remark}
\newtheorem{remarks}[proposition]{Remarks}
\newcommand{\thlabel}[1]{\label{th:#1}}
\newcommand{\thref}[1]{Theorem~\ref{th:#1}}
\newcommand{\selabel}[1]{\label{se:#1}}
\newcommand{\seref}[1]{Section~\ref{se:#1}}
\newcommand{\lelabel}[1]{\label{le:#1}}
\newcommand{\leref}[1]{Lemma~\ref{le:#1}}
\newcommand{\prlabel}[1]{\label{pr:#1}}
\newcommand{\prref}[1]{Proposition~\ref{pr:#1}}
\newcommand{\colabel}[1]{\label{co:#1}}

\newcommand{\relabel}[1]{\label{re:#1}}

\newcommand{\reslabel}[1]{\label{res:#1}}
\newcommand{\resref}[1]{Remark~\ref{res:#1}}

\newcommand{\delabel}[1]{\label{de:#1}}
\newcommand{\deref}[1]{Definition~\ref{de:#1}}
\newcommand{\eqlabel}[1]{\label{eq:#1}}
\newcommand{\equref}[1]{(\ref{eq:#1})}
\def\ra{\rightarrow}
\def\cd{\cdot}
\newcommand{\btrl}
{\mbox{${\;}$$\blacktriangleright\hspace{-0.8mm}\blacktriangleleft$${\;}$}}
\newcommand\smi{\mbox{$S^{-1}$}}
\def\ot{\otimes}
\def\va{\varepsilon}

\def\un{\underline}
\def\mb{\mathbb}
\def\mf{\mathfrak}

\def\mfA{\mf {A}}
\def\mfa{\mf {a}}
\def\le{\langle}
\def\ri{\rangle}
\def\l{\lambda}
\def\r{\rho}

\def\va{\varepsilon}
\def\v{\varphi}

\def\rh{\rightharpoonup}
\def\lh{\leftharpoonup}

\def\ra{\rightarrow}
\def\a{\alpha}
\def\b{\beta}

\def\d{\delta}
\def\O{\Omega}

\def\ov{\overline}
\def\cal{\mathcal}
\def\un{\underline}

\newcommand{\mfb}{\mbox{$\mf {b}$}}

\newcommand{\tx}{\mbox{$\tilde {x}$}}
\newcommand{\tX}{\mbox{$\tilde {X}$}}
\newcommand{\ty}{\mbox{$\tilde {y}$}}

\newcommand{\tpla}{\mbox{$\tilde {p}^1_{\l }$}}
\newcommand{\tplb}{\mbox{$\tilde {p}^2_{\l }$}}
\newcommand{\tPla}{\mbox{$\tilde {P}^1_{\l }$}}
\newcommand{\tPlb}{\mbox{$\tilde {P}^2_{\l }$}}
\newcommand{\tqla}{\mbox{$\tilde {q}^1_{\l }$}}
\newcommand{\tqlb}{\mbox{$\tilde {q}^2_{\l }$}}
\newcommand{\tQla}{\mbox{$\tilde {Q}^1_{\l }$}}
\newcommand{\tQlb}{\mbox{$\tilde {Q}^2_{\l }$}}

\newcommand{\und}{\mbox{$\un {\Delta }$}}
\newcommand{\une}{\mbox{$\un {\va }$}}
\newcommand{\una}{\mbox{$c_{\un {1}}$}}
\newcommand{\unb}{\mbox{$c_{\un {2}}$}}
\newcommand{\unaa}{\mbox{$c_{(\un {1}, \un {1})}$}}
\newcommand{\unab}{\mbox{$c_{(\un {1}, \un {2})}$}}

\newcommand{\Hom}{{\rm Hom}}

\newcommand{\im}{{\rm Im}\,}

\newcommand{\Rat}{{\rm Rat}}

\def\lan{\langle}
\def\ran{\rangle}
\def\ot{\otimes}

\def\cal{\mathcal}

\newcommand{\Cc}{\mathcal{C}}

\newcommand{\Mm}{\mathcal{M}}

\def\*C{{}^*\hspace*{-1pt}{\Cc}}
\def\text#1{{\rm {\rm #1}}}

\def\ul{\underline}

\newcommand{\gsm}{\mbox{$\blacktriangleright \hspace{-0.7mm}<$}}
\newcommand{\gtl}{\mbox{${\;}$$>\hspace{-0.85mm}\blacktriangleleft$${\;}$}}
\begin{document}
\title[Doi-Hopf modules and Yetter-Drinfeld modules]
{Doi-Hopf modules and Yetter-Drinfeld modules for quasi-Hopf algebras}
\author{D. Bulacu}
\address{Faculty of Mathematics and Informatics, University 
of Bucharest,RO-010014 Bucha-rest 1, Romania}
\email{dbulacu@al.math.unibuc.ro}
\thanks{ 
The first author was financially supported by the bilateral project
"Hopf Algebras in Algebra, Topology, Geometry and Physics" of the 
Flemish and Romanian governments, and SB 2002-0286. He would like to 
thank Free University of Brussels (Belgium) and University of Almeria 
(Spain) for their warm hospitality.}
\author{S. Caenepeel}
\address{Faculty of Applied Sciences, 
Vrije Universiteit Brussel, VUB, B-1050 Brussels, Belgium}
\email{scaenepe@vub.ac.be}
\urladdr{http://homepages.vub.ac.be/\~{}scaenepe/}
\author{B. Torrecillas}
\address{Department of Algebra and Analysis, 
University of Almeria, 04071 Almeria, Spain}
\email{btorreci@ual.es}
\subjclass{16W30}
\keywords{quasi-Hopf algebra, Doi-Hopf module, Yetter-Drinfeld module}

\begin{abstract}
For a quasi-Hopf algebra $H$, a left $H$-comodule algebra 
$\mf{B}$ and a right $H$-module coalgebra $C$ we will characterize the 
category of Doi-Hopf modules ${}^C{\cal M}(H)_{\mf{B}}$ in terms 
of modules. We will also show that for an 
$H$-bicomodule algebra $\mb{A}$ and an $H$-bimodule coalgebra 
$C$ the category of generalized Yetter-Drinfeld modules 
${}_{\mb{A}}{\cal YD}(H)^C$ is isomorphic to a 
certain category of Doi-Hopf modules. Using 
this isomorphism we will transport the properties from the category of 
Doi-Hopf modules to the category of generalized Yetter-Drinfeld modules.  
\end{abstract}
\maketitle
%%%%%%%%%%%%%%%%%%%%%%%%%%%%%%%
\section*{Introduction}
%%%%%%%%%%%%%%%%%%%%%%%%%%%%%%%%%%
Recall that the defining axioms for a quasi-bialgebra $H$ are the same as
for a bialgebra, with the coassociativity of the comultiplication
replaced by a weaker property, called quasi-coassociativity: the comultiplication
is coassociative up to conjugation by an invertible element
$\Phi \in H\ot H\ot H$, called the reassociator. There are important differences
with ordinary quasi-bialgebras: the definition of a quasi-bialgebra is
not selfdual, and we cannot consider comodules over quasi-bialgebras,
since they are not coassociative coalgebras. However, the category of
(left or right) modules over a quasi-bialgebra is a monoidal category.

Using this categorical point of view, the category of relative 
Hopf modules has been introduced and studied in \cite{bn}. A right
$H$-module coalgebra $C$ is a coalgebra in the monoidal category
${\cal M}_H$, and a left $[C, H]$-Hopf module is a left $C$-comodule 
in the monoidal category ${\cal M}_H$. 
A generalization of this concept was presented in \cite{bc3}:
replacing the right $H$-action by an action of a left $H$-comodule algebra, 
we can define the notion of Doi-Hopf module over a quasi-bialgebra.
At this point, we have to mention that there is a phylosophical problem
with the introduction of $H$-comodule algebras: we cannot introduce them
as algebras in the category of $H$-comodules, since this category does not
exist, as we mentioned above. However, a formal definition of $H$-comodule
(and $H$-bicomodule)
algebra was given by Hausser and Nill in \cite{hn1}. A more conceptual
definition 
has been proposed in \cite{bc3}. If $\mf{A}$ is an associative algebra 
then the category of $(\mf{A}\ot H, \mf{A})$-bimodules is monoidal. Moreover, 
$\mf{A}\ot H$ has a coalgebra structure within this monoidal category 
``compatible" with the unit element $1_{\mf{A}}\ot 1_H$ if and only if 
$\mf{A}$ is a right $H$-comodule algebra (for the complete statement see 
\prref{1.2} below). Of course, a similar result holds for a left $H$-comodule algebra 
(see \prref{1.3}). Moreover, if $H$ is a quasi-Hopf algebra then 
any $H$-bicomodule algebra can be viewed in two different (but twist equivalent) 
ways as a left (right) $H$-comodule algebra.

The aim of this paper is to study the category of Doi-Hopf modules over a 
quasi-Hopf algebra $H$, and its connections to the category of 
Yetter-Drinfeld modules. If $H$ is a quasi-bialgebra, 
$\mf{B}$ a left $H$-comodule algebra 
and $C$ a right $H$-module coalgebra then $\mf{B}\ot C$ is a 
$\mf{B}$-coring (this means a coalgebra in the monoidal category 
of $\mf{B}$-bimodules) and the category of right-left 
$(H, \mf{B}, C)$-Hopf modules, denoted by ${}^C{\cal M}(H)_{\mf{B}}$, 
is isomorphic to the category of right comodules over the coring $\mf{B}\ot C$, 
conform \cite[Theorem 5.4]{bc3}. In particular, since $\mf{B}\ot C$ is flat as a 
left $\mf{B}$-module we obtain 
that ${}^C{\cal M}(H)_{\mf{B}}$ is a Grothendieck category. 

It was shown in \cite[Proposition 5.2]{bc3} that the category
${}^C{\cal M}(H)_{\mf{B}}$ is isomorphic to the category of right modules 
over the generalized smash product $C^*\gsm \mf{B}$, if 
$C$ is finite dimensional. 
In \seref{2} we look at the case where $C$ is infinite dimensional. Following the methods 
developed in \cite{dnt,doi,mz}, we will present two characterizations
of the category of Doi-Hopf modules  ${}^C{\cal M}(H)_{\mf{B}}$. 
We first introduce the notion 
of rational (right) $C^*\gsm \mf{B}$-module, and then we will show that 
the category 
${}^C{\cal M}(H)_{\mf{B}}$ is isomorphic to $\Rat({\cal M}_{C^*\gsm \mf{B}})$, 
the category of rational (right) $C^*\gsm \mf{B}$-modules. We notice that, in the 
coassociative case, the notion of rational (right) $C^*\gsm \mf{B}$-module reduces 
to the notion of right $C^*\gsm \mf{B}$-module which is rational as a 
right $C^*$-module. 
Secondly, we will show that 
%${}^C{\cal M}(H)_{\mf{B}}$ is equivalent to 
$\Rat({\cal M}_{C^*\gsm \mf{B}})$ is equal to 
$\sigma _{C^*\gsm \mf{B}}[C\ot \mf{B}]$, the smallest closed subcategory of 
${\cal M}_{C^*\gsm \mf{B}}$ containing $C\ot \mf{B}$ (see \thref{2.5}). 
In this way we recover that ${}^C{\cal M}(H)_{\mf{B}}$ is a Grothendieck category, 
a fortiori with enough injective objects. We will also introduce a generalized
version of Koppinen's smash product \cite{koppinen}, relate it to the
generalized smash product, and characterize the category of Doi-Hopf modules
as the full subcategory of modules over the Koppinen smash product,
consisting of rational modules.

In \seref{3}, we will generalize a result from \cite{cmz}. 
If $H$ is a quasi-Hopf algebra, then an $H$-bicomodule algebra 
$\mb{A}$ can be viewed in two different, but twist equivalent, ways as 
a right $H^{\rm op}\ot H$-comodule algebra. To this end, we first 
prove that any left $H$-comodule algebra $\mf{B}$ can be turned into 
a right $H^{\rm op}$-comodule algebra. So, by this correspondence, the two 
(twist equivalent) left $H\ot H^{\rm op}$-comodule algebra structures 
on $\mb{A}$ obtained in \cite{bpv3} provide two (twist equivalent) 
right $H^{\rm op}\ot H$-comodule algebra structures on $\mb{A}$, 
which we will denote by $\mb{A}^1$ and $\mb{A}^2$. If $C$ is an $H$-bimodule 
coalgebra (that is, a coalgebra in the monoidal category of $H$-bimodules), 
then $C$ becomes in a natural way a left $H^{\rm op}\ot H$-module coalgebra, thus it 
makes sense to consider the Hopf module category 
${}_{\mb{A}^2}{\cal M}(H^{\rm op}\ot H)^C$. The main result of \seref{3} asserts 
that the category of generalized left-right Yetter-Drinfeld modules 
${}_{\mb{A}}{\cal YD}(H)^C$ is isomorphic to 
${}_{\mb{A}^2}{\cal M}(H^{\rm op}\ot H)^C$, and also to 
${}_{\mb{A}^1}{\cal M}(H^{\rm op}\ot H)^C$. Using the 
first isomorphism, we will characterize ${}_{\mb{A}}{\cal YD}(H)^C$ 
as a category of comodules over a coring. In \seref{3.3}, we will characterize
the category of Yetter-Drinfeld modules as a category of modules.

%%%%%%%%%%%%%%%%%%%%%%%%%%%%%%%
\section{Preliminary results}\selabel{1}
%%%%%%%%%%%%%%%%%%%%%%%%%%%%%%%%%
\subsection{Quasi-Hopf algebras}\selabel{1.1}
%%%%%%%%%%%%%%%%%%%%%%%%%%%%%%%%%
We work over a commutative field $k$. All algebras, linear spaces
etc. will be over $k$; unadorned $\ot $ means $\ot_k$. Following
Drinfeld \cite{d1}, a quasi-bialgebra is a fourtuple 
$(H, \Delta , \va , \Phi )$, where $H$ is an associative algebra with unit, 
$\Phi$ is an invertible element in $H\ot H\ot H$, and 
$\Delta :\ H\ra H\ot H$ and $\va :\ H\ra k$ are algebra homomorphisms
satisfying the identities 
\begin{eqnarray}
&&(id \ot \Delta )(\Delta (h))=
\Phi (\Delta \ot id)(\Delta (h))\Phi ^{-1},\label{q1}\\
&&(id \ot \va )(\Delta (h))=h, ~~ 
(\va \ot id)(\Delta (h))=h,\label{q2}
\end{eqnarray}
for all $h\in H$; $\Phi$ has to be a normalized $3$-cocycle,
in the sense that
\begin{eqnarray}
&&(1\ot \Phi)(id\ot \Delta \ot id) (\Phi)(\Phi \ot 1)
=(id\ot id \ot \Delta )(\Phi ) (\Delta \ot id \ot id)
(\Phi ),\label{q3}\\
&&(id \ot \va \ot id )(\Phi )=1\ot 1.\label{q4}
\end{eqnarray} 
The map $\Delta $ is called the coproduct or the
comultiplication, $\va $ the counit and $\Phi $ the reassociator.
As for Hopf algebras \cite{sw} we denote $\Delta (h)=h_1\ot h_2$, 
but since $\Delta $ is only quasi-coassociative we adopt the 
further convention (summation implicitely understood):
$$
(\Delta \ot id)(\Delta (h))=h_{(1, 1)}\ot h_{(1, 2)}\ot h_2, ~~ 
(id\ot \Delta )(\Delta (h))=h_1\ot h_{(2, 1)}\ot h_{(2,2)}, 
$$
for all $h\in H$. We will denote the tensor components of $\Phi$
by capital letters, and those of $\Phi^{-1}$ by small letters, 
namely
\begin{eqnarray*}
&&\Phi=X^1\ot X^2\ot X^3=T^1\ot T^2\ot T^3
=V^1\ot V^2\ot V^3=\cdots\\
&&\Phi^{-1}=x^1\ot x^2\ot x^3=t^1\ot t^2\ot t^3
=v^1\ot v^2\ot v^3=\cdots
\end{eqnarray*}

A quasi-Hopf algebra is a quasi-bialgebra $H$ equipped with an 
anti-automorphism $S$ of the algebra $H$ and elements $\a , \b \in
H$ such that
\begin{eqnarray}
&&S(h_1)\a h_2=\va (h)\a {~~\rm and ~~}
h_1\b S(h_2)=\va (h)\b ,\label{q5}\\[1mm]
&&X^1\b S(X^2)\a X^3=1 {~~\rm and ~~}
S(x^1)\a x^2\b S(x^3)=1.\label{q6},
\end{eqnarray}
for all $h\in H$.

The antipode of a quasi-Hopf algebra is determined
uniquely up to a transformation $\a \mapsto U\a $, $\b \mapsto \b
U^{-1}$, $S(h)\mapsto US(h)U^{-1}$, with $U\in H$ invertible.
The axioms imply that 
$\va (\a )\va (\b )=1$, so, by rescaling  
$\a $ and $\b $, we may assume without loss 
of generality that $\va (\a )=\va (\b )=1$ and $\va \circ S=\va $.
The identities (\ref{q2}-\ref{q4}) also imply 
that
$$
(\va \ot id\ot id)(\Phi )=(id \ot id\ot \va )(\Phi )=1\ot 1.
$$

If $H=(H, \Delta , \va , \Phi , S, \a , \b )$ is 
a quasi-bialgebra or a quasi-Hopf algebra then 
$H^{\rm op}$, $H^{\rm cop}$ and $H^{\rm op,cop}$ are also quasi-bialgebras 
(respectively quasi-Hopf algebras), 
where the superscript ``op" means opposite 
multiplication and ``cop" means opposite comultiplication. The 
structure maps are obtained by putting $\Phi _{\rm op}=\Phi ^{-1}$,
$\Phi _{\rm cop}=(\Phi ^{-1})^{321}$, $\Phi _{\rm op,cop}=\Phi ^{321}$,
$S_{\rm op}=S_{\rm cop}=(S_{\rm op,cop})^{-1}=S^{-1}$, $\a _{\rm op}=\smi (\b )$,
$\b _{\rm op}=\smi (\a )$, $\a _{\rm cop}=\smi (\a )$, $\b _{\rm cop}=\smi (\b )$,
$\a _{\rm op,cop}=\b $ and $\b _{\rm op,cop}=\a $.

The definition of a quasi-bialgebra $H$ is designed in such a way that 
the categories of left and right representations over $H$ are monoidal
(see \cite{k, m2} for the terminology). Let $(H, \Delta , \varepsilon , \Phi )$ be a
quasi-bialgebra. For  are left ( resp. right) $H$-modules $U,V,W$, 
the associativity constraints 
$a_{U,V,W}~ ({\rm resp.~} {\bf a}_{U, V, W}) :\ (U\otimes V)\otimes W\rightarrow
U\otimes (V\otimes W)$
are given by the formulas
$$a_{U,V,W}((u\otimes v)\otimes w)=\Phi \cdot (u\otimes
(v\otimes w));$$
$${\bf a}_{U, V, W}((u\ot v)\ot w)= (u\ot (v\ot w))\cd \Phi ^{-1}.$$

Next we recall that the definition of a quasi-bialgebra or 
quasi-Hopf algebra is 
``twist covariant" in the following sense. An invertible element
$F\in H\ot H$ is called a gauge transformation or 
twist if $(\va \ot id)(F)=(id\ot \va)(F)=1$. If $H$ is a quasi-bialgebra 
or a quasi-Hopf algebra and $F=F^1\ot F^2\in H\ot H$ is a gauge 
transformation with inverse $F^{-1}=G^1\ot G^2$, then we can
define a new quasi-bialgebra (respectively quasi-Hopf algebra) 
$H_F$ by keeping the 
multiplication, unit, counit (and antipode in the case of a quasi-Hopf 
algebra) of $H$ and replacing the 
comultiplication, reassociator and the elements $\alpha$ and $\beta$ by 
\begin{eqnarray}
&&\Delta _F(h)=F\Delta (h)F^{-1},\label{g1}\\[1mm]
&&\Phi_F=(1\ot F)(id \ot \Delta )(F) \Phi (\Delta \ot id)
(F^{-1})(F^{-1}\ot 1),\label{g2}\\[1mm]
&&\a_F=S(G^1)\a G^2, ~~ 
\b_F=F^1\b S(F^2).\label{g3}
\end{eqnarray}

It is well-known that the antipode of a Hopf
algebra is an anti-coalgebra morphism. For a quasi-Hopf algebra,
we have the following statement: there exists a gauge
transformation $f\in H\ot H$ such that
\begin{equation} \label{ca}
f\Delta (S(h))f^{-1}=(S\ot S)(\Delta ^{\rm cop}(h)), ~~
\mbox{for all $h\in H$.}
\end{equation}
The element $f$ can be computed explicitly. First set 
\begin{eqnarray*}
&&A^1\ot A^2\ot A^3\ot A^4=(\Phi \ot 1) (\Delta \ot id\ot id)(\Phi ^{-1}),\\
&&B^1\ot B^2\ot B^3\ot B^4=(\Delta \ot id\ot id)(\Phi )(\Phi ^{-1}\ot 1),
\end{eqnarray*}
and then define $\gamma, \delta\in H\ot H$ by
$$
\gamma =S(A^2)\a A^3\ot S(A^1)\a A^4~~{\rm and}~~ 
\delta =B^1\b S(B^4)\ot B^2\b S(B^3).
$$
Then $f$ and $f^{-1}$ are given by the formulae
\begin{eqnarray}
f&=&(S\ot S)(\Delta ^{\rm cop}(x^1)) \gamma \Delta (x^2\b
S(x^3)),\label{f}\\%
f^{-1}&=&\Delta (S(x^1)\a x^2) \delta (S\ot S)(\Delta
^{\rm cop}(x^3)).\label{g}
\end{eqnarray}
Furthermore the corresponding twisted reassociator (see
(\ref{g2})) is given by
\begin{equation} \label{pf}
\Phi _f=(S\ot S\ot S)(X^3\ot X^2\ot X^1).
\end{equation}

%%%%%%%%%%%%%%%%%%%%%%%%%%%%%%%%%%%%%%
\subsection{Comodule and bicomodule algebras}\selabel{1.2}
%%%%%%%%%%%%%%%%%%%%%%%%%%%%%%%%%%%%%%%%%%%%%
A formal definition of comodule algebras over a 
quasi-bialgebra was given by Hausser and Nill \cite{hn1}.

\begin{definition}\delabel{1.1}
Let $H$ be a quasi-bialgebra. A right $H$-comodule algebra
is a unital associative algebra
$\mathfrak{A}$ together with an algebra morphism $\r :\mathfrak{A}\ra \mathfrak{A}\ot H$
and an invertible element $\Phi _{\r }\in \mathfrak{A}\ot H\ot H$
such that:
\begin{eqnarray}
&&\hspace*{-2cm}\Phi _{\r }(\r \ot id)(\r (\mf {a}))
=(id\ot \Delta )(\r (\mf {a}))\Phi _{\r }, 
~~{\rm for~all} ~\mf {a}\in \mathfrak{A},\label{rca1}\\
&&\hspace*{-2cm}(1_{\mf {A}}\ot \Phi)(id\ot \Delta \ot id)(\Phi _{\r })
(\Phi _{\r }\ot 1_H)\nonumber\\
&=&(id\ot id\ot \Delta )(\Phi _{\r })
(\r \ot id\ot id)(\Phi _{\r }),\label{rca2}\\
&&\hspace*{-2cm}(id\ot \va)\circ \r =id ,\label{rca3}\\
&&\hspace*{-2cm}(id\ot \va \ot id)(\Phi _{\r })=(id\ot id\ot \va )
(\Phi _{\r })=1_{\mathfrak{A}}\ot 1_H.\label{rca4}
\end{eqnarray}

In a similar way, a left $H$-comodule algebra is a
unital associative algebra $\mathfrak{B}$ together with
an algebra morphism 
$\l : \mf {B}\ra H\ot \mathfrak{B}$ and an invertible element 
$\Phi _{\l }\in H\ot H\ot \mathfrak{B}$ such that the following 
relations hold:
\begin{eqnarray}
&&\hspace*{-2cm}(id\ot \l )(\l (\mf {b}))\Phi _{\l }=
\Phi _{\l }(\Delta \ot id)(\l (\mf {b})),
~~\forall ~\mf {b}\in \mathfrak{B},\label{lca1}\\
&&\hspace*{-2cm}(1_H\ot \Phi _{\l })(id\ot \Delta \ot id)(\Phi _{\l })
(\Phi \ot 1_{\mf {B}})\nonumber\\
&=&(id\ot id\ot \l )(\Phi _{\l })
(\Delta \ot id\ot id)(\Phi _{\l }),\label{lca2}\\
&&\hspace*{-2cm}(\va \ot id)\circ \l =id ,\label{lca3}\\
&&\hspace*{-2cm}(id\ot \va \ot id)(\Phi _{\l })=(\va \ot id\ot id)
(\Phi _{\l })=1_H\ot 1_{\mathfrak{B}}.\label{lca4}
\end{eqnarray}
\end{definition}

Observe that $H$ is a left and right $H$-comodule algebra: take
$\r =\l =\Delta $, 
$\Phi _{\r }=\Phi _{\l }=\Phi $. 

If $(\mf{B}, \l , \Phi _{\l})$ is a left $H$-comodule algebra then 
\begin{itemize}
\item[-] $(\mf{B}, \l \circ \tau _{H, \mf{B}}, (\Phi _{\l }^{-1})^{321})$ 
is a right $H^{\rm cop}$-comodule algebra;  
\item[-] $(\mf {B}^{\rm op}, \l \circ \tau _{H, \mf{B}}, \Phi _{\l }^{321})$ is a 
right $H^{\rm op,cop}$-comodule algebra; 
\item[-] $(\mf {B}^{\rm op}, \l , \Phi _{\l }^{-1})$ is a left 
$H^{\rm op}$-comodule algebra, 
\end{itemize}
and vice versa. $\tau_{X,Y}:\ X\ot Y\to Y\ot X$ is the switch map mapping
$x\ot y$ to $y\ot x$.

From \cite{hn1}, we recall the notion of twist equivalence for the
coaction on a right $H$-comodule algebra $(\mfA , \r , \Phi _{\r})$:
if $\mb{V}\in \mfA \ot H$ is an invertible element such that 
\begin{equation}\eqlabel{comtwist0}
(id_{\mfA}\ot \va )(\mb{V})=1_{\mfA}
\end{equation}
 then we can construct a
new right $H$-comodule algebra $(\mfA , \r ', \Phi _{\r '})$ with
\begin{equation}\eqlabel{comtwist1}
\r '(\mfa )=\mb{V}\r (\mfa )\mb{V}^{-1}
\end{equation}
and
\begin{equation}\eqlabel{comtwist2}
\Phi _{\r '}=(id_{\mfA}\ot \Delta )(\mb{V})\Phi _{\r }(\r \ot id_H)
(\mb{V}^{-1})(\mb{V}^{-1}\ot 1_H).
\end{equation}
We say that $(\mfA , \r , \Phi _{\r})$ and $(\mfA , \r ', \Phi _{\r}')$ 
are twist equivalent right $H$-comodule algebras. 

We obviously have a similar notion for left $H$-comodule algebras.
More precisely, if 
$(\mf{B}, \l , \Phi _{\l})$ is a left $H$-comodule algebra and 
$\mb{U}\in H\ot \mf{B}$ is an invertible element such that 
$(\va \ot id_{\mf{B}})(\mb{U})=1_{\mf{B}}$ then we have a new
left $H$-comodule algebra
$(\mf{B}, \l ', \Phi _{\l '})$ with
$\l '(\mfb )=\mb{U}\l (\mfb )\mb{U}^{-1}$, for all $\mfb \in \mf{B}$, and
$$\Phi _{\l '}=(1_H\ot \mb{U})(id_H\ot \l )(\mb{U})\Phi _{\l }
(\Delta \ot id_{\mf{B}})\mb{U}^{-1}.$$
$(\mf{B}, \l , \Phi _{\l})$ 
and $(\mf{B}, \l ', \Phi _{\l}')$ are called twist equivalent. 

For a right $H$-comodule algebra $({\mf A}, \r , \Phi _{\r })$ we
will use the following Sweedler-type notation, for any $\mfa \in {\mf A}$. 
$$
\r (\mfa )=\mfa _{\le0\ri}\ot \mfa _{\le1\ri}, ~~
(\r \ot id)(\r (\mfa ))=\mfa _{\le0, 0\ri}\ot \mfa _{\le0, 1\ri}\ot 
\mfa _{\le1\ri} \mbox{~~etc.}
$$
Similarly, for a left $H$-comodule 
algebra $({\mf B}, \l , \Phi _{\l })$, if $\mfb \in {\mf B}$, we adopt 
the following notation: 
$$
\l (\mfb )=\mfb _{[-1]}\ot \mfb _{[0]}, ~~ 
(id\ot \l )(\l (\mfb ))=\mfb _{[-1]}\ot \mfb _{[0,-1]}\ot
\mfb _{[0, 0]} \mbox{~~etc.} 
$$
In analogy with the notation for the reassociator $\Phi $ of $H$, 
we will write 
$$
\Phi _{\r }=\tilde {X}^1_{\r }\ot \tilde {X}^2_{\r }
\ot \tilde {X}^3_{\r }=
\tilde {Y}^1_{\r }\ot \tilde {Y}^2_{\r }\ot 
\tilde {Y}^3_{\r }=\cdots  
$$ 
and 
$$ 
\Phi _{\r }^{-1}=\tilde {x}^1_{\r }\ot \tilde {x}^2_{\r }\ot
\tilde {x}^3_{\r }=\tilde {y}^1_{\r }\ot 
\tilde {y}^2_{\r }\ot \tilde {y}^3_{\r }=\cdots .
$$ 
We use a similar notation for the element $\Phi _{\l }$ of a left $H$-comodule 
algebra $\mf {B}$. 

If $H$ is an associative bialgebra and $\mf {A}$ is an ordinary right 
$H$-comodule algebra, then $\mf {A}\ot H$ is a coalgebra 
in the monoidal category of $\mf {A}$-bimodules. The 
quasi-bialgebra analog of this property was given in \cite{bc3}. 
Let $H$ be a quasi-bialgebra and $\mf {A}$ a unital associative 
algebra. We define by ${}_{\mf {A}\ot H}{\cal M}_{\mf {A}}$ 
the category of $\mf {A}$-bimodules and $(H, \mf {A})$-bimodules 
$M$ such that $h\cd (\mfa \cd m)=\mfa \cd (h\cd m)$, for all $\mfa \in \mf {A}$, 
$h\in H$ and $m\in M$. Morphisms are left $H$-linear maps which are
also $\mf {A}$-bimodule maps. It is not hard to see that 
${}_{\mf {A}\ot H}{\cal M}_{\mf {A}}$ is a monoidal category. 
The tensor product is $\ot _{\mf {A}}$ and for any two objects 
$M, N\in {}_{\mf {A}\ot H}{\cal M}_{\mf {A}}$, $M\ot _{\mf {A}}N$ 
is an object of ${}_{\mf {A}\ot H}{\cal M}_{\mf {A}}$ via 
$$
(\mfa \ot h)\cd (m\ot _{\mf {A}}n)\cd \mfa {'}=
\mfa \cd (h_1\cd m)\ot _{\mf {A}}h_2\cd n\cd \mfa {'} 
$$
for all $m\in M$, $n\in N$, $\mfa , \mfa {'}\in \mf {A}$, 
and $h\in H$. The associativity constraints are given by 
\begin{eqnarray*}
&&\un{a}_{M, N, P}: (M\ot _{\mf A}N)
\ot _{\mf A}P\ra M\ot _{\mf A}(N\ot _{\mf A}P),\\
&&\un{a}_{M, N, P}((m\ot_{\mf A}n)\ot_{\mf A}p)=
X^1\cd m\ot _{\mf A}(X^2\cd n\ot _{\mf A}X^3\cd p);
\end{eqnarray*}
the unit object is $\mf A$ viewed as a trivial left $H$-module, and the 
left and right unit constraints are the usual ones. Now, the definition 
of a comodule algebra in terms of monoidal categories can be 
restated as follows.

\begin{proposition}{\rm (\cite[Proposition 3.8]{bc3})}\prlabel{1.2}
Let $H$ be a quasi-bialgebra and $\mf A$ an algebra, and view
${\mf A}\ot H$ in the canonical way as an
object in ${}_{{\mf A}\ot H}{\cal M}$. There is a bijective correspondence
between coalgebra structures $({\mf A}\ot H, \und , \une )$
in the monoidal category ${}_{{\mf A}\ot H}{\cal M}_{\mf A}$
such that $\und (1_{\mf A}\ot 1_H)$ is invertible and
$\une (1_{\mf A}\ot 1_H)=1_{\mf A}$, and  right
$H$-comodule algebra structures on $\mf A$.
\end{proposition}

A similar result holds for a left $H$-comodule algebra $\mf{B}$. 
Denote by ${}_{\mf{B}}{\cal M}_{\mf{B}\ot H}$ the category whose objects are 
$\mf{B}$-bimodules and $(\mf{B}, H)$-bimodules $M$ such that 
$(m\cd \mfb )\cd h=(m\cd h)\cd \mfb $ for all $m\in M$, $\mfb \in \mf{B}$ 
and $h\in H$. Morphisms are right $H$-linear maps which are also $\mf{B}$-bimodule 
maps. We can easily check that ${}_{\mf{B}}{\cal M}_{\mf{B}\ot H}$ is a monoidal 
category with tensor product $\ot _{\mf{B}}$ given via $\Delta $, this means 
$$
\mfb \cd (m\ot _{\mf{B}}n)\cd (\mfb '\ot h):=\mfb \cd m\cd h_1
\ot _{\mf{B}}(n\cd h_2)\cd \mfb '
$$    
for all $M, N\in {}_{\mf{B}}{\cal M}_{\mf{B}\ot H}$, $m\in M$, $n\in N$, 
$\mfb , \mfb '\in \mf{B}$ and $h\in H$. The associativity constraints are 
given by
\begin{eqnarray*}
&&\un{a}'_{M, N, P}:\ (M\ot _{\mf{B}}N)\ot _{\mf{B}}P\ra 
M\ot _{\mf{B}}(N\ot _{\mf{B}}P),\\
&&\un{a}'_{M, N, P}((m\ot _{\mf{B}}n)\ot _{\mf{B}}p)=
m\cd x^1\ot _{\mf{B}}(n\cd x^2\ot _{\mf{B}}p\cd x^3),
\end{eqnarray*}
the unit object is $\mf{B}$ viewed as a right $H$-module via $\va $, and 
the left and right unit constraints are the usual ones. 

\begin{proposition}\prlabel{1.3}
Let $H$ be a quasi-bialgebra and $\mf B$ an algebra, and view
${\mf B}\ot H$ in the canonical way as an
object in ${\cal M}{}_{{\mf B}\ot H}$. There is a bijective correspondence
between coalgebra structures $(\mf{B}\ot H, \und ,\une )$
in the monoidal category ${}_{\mf{B}}{\cal M}_{\mf{B}\ot H}$ 
such that $\und (1_{\mf{B}}\ot 1_H)$ 
is invertible and $\une (1_{\mf{B}}\ot 1_H)=1_{\mf{B}}$, and left
$H$-comodule algebra structures on $\mf B$.
\end{proposition} 

\begin{proof}
Since the proof is similar to the one given in \cite[Proposition 3.8]{bc3},
we restrict ourselves to a description of the correspondence.

Suppose that $(\mf{B}\ot H, \und ,\une )$ is a coalgebra in
${}_{\mf{B}}{\cal M}_{\mf{B}\ot H}$ such that $\und (1_{\mf{B}}\ot 1_H)$ 
is invertible and $\une (1_{\mf{B}}\ot 1_H)=1_{\mf{B}}$. Write 
$$
\und (1_{\mf{B}}\ot 1_H)=(1_{\mf{B}}\ot \tilde{X}^1_{\l})
\ot _{\mf{B}}(\tilde{X}^3_{\l}\ot \tilde{X}^2_{\l}),
$$
and consider 
$\Phi _{\l}=\tilde{X}^1_{\l}\ot \tilde{X}^2_{\l}\ot \tilde{X}^3_{\l}$. Define
$\lambda:\ \mf\to H\ot \mf{B}$ by
$$\lambda(b)=\tau _{\mf{B}, H}(\mfb \cd (1_{\mf{B}}\ot 1_H)),$$
for all $\mfb \in \mf{B}$. Then $(\mf{B}, \l , \Phi _{\l })$ is a 
left $H$-comodule algebra. 

Conversely, if $(\mf{B}, \l , \Phi _{\l})$ is a left $H$-comodule algebra then 
$\mf{B}\ot H\in {}_{\mf{B}}{\cal M}_{\mf{B}\ot H}$ via
$$
\mfb \cd (\mfb '\ot h)\cd (\mfb ''\ot h')=\mfb _{[0]}\mfb '\mfb ''\ot 
\mfb _{[-1]}hh'.
$$
Moreover, $\mf{B}\ot H$ is a coalgebra in 
${}_{\mf{B}}{\cal M}_{\mf{B}\ot H}$, with comultiplication and counit
given by the formulas
$$\und (\mfb \ot h)=(1_{\mf{B}}\ot \tilde{X}^1_{\l}h_1)\ot _{\mf{B}}
(\tilde{X}^3_{\l}\mfb \ot \tilde{X}^2_{\l}h_2),$$
$$\une (\mfb \ot h)=\va (h)\mfb,$$
for all $\mfb \in \mf{B}$ and $h\in H$.
\end{proof}

Let $\mf {B}$ be a left $H$-comodule algebra, and consider
the elements $\tilde{p}_{\l}$ and 
$\tilde{q}_{\l }$ in $H\ot \mf {B}$ given by the following formulas:
\begin{equation}\label{tpql}
\tilde{p}_{\l }=\tpla \ot \tplb 
=\tilde {X}^2_{\l }\smi (\tilde{X}^1_{\l }\b )\ot 
\tilde{X}^3_{\l }, ~~
\tilde{q}_{\l }=S(\tilde{x}^1_{\l })\a \tilde{x}^2_{\l }
\ot \tilde{x}^3_{\l }.
\end{equation}
Then we have the following formulas, for all $\mfb \in \mf {B}$ (see \cite{hn1}):
\begin{eqnarray}
&&\l (\mfb _{[0]})\tilde{p}_{\l }[\smi (\mfb _{[-1]})\ot 1_{\mf {B}}]
=\tilde{p}_{\l }[1_H\ot \mfb ],\label{tpql1}\\
&&[S(\mfb _{[-1]})\ot 1_H]\tilde{q}_{\l }\l (\mfb _{[0]})=
[1\ot \mfb ]\tilde{q}_{\l },\label{tpql1a}\\
&&\l (\tqlb )\tilde{p}_{\l }[\smi (\tqla )\ot 1_{\mf {B}}]
=1_H\ot 1_{\mf {B}},\label{tpql2}\\
&&[S(\tpla )\ot 1_{\mf {B}}]\tilde{q}_{\l }\l (\tplb )
=1_H\ot 1_{\mf {B}},\label{tpql2a}\\
&&\Phi _{\l }^{-1}(id_{\mf{B}}\ot \l)(\tilde{p}_{\l})
(1_{\mf{B}}\ot \tilde{p}_{\l})\nonumber\\ 
&&\hspace*{1cm}=
(\Delta \ot id_{\mf{B}})(\l (\tilde{X}^3_{\l})\tilde{p}_{\l})
[\smi (\tilde{X}^2_{\l}g^2)\ot 
\smi (\tilde{X}^1_{\l}g^1)\ot 1_{\mf{B}}],\label{tpl}\\
&&(1_H\ot \tilde{q}_{\l})(id_H\ot \l )(\tilde{q}_{\l})\Phi _{\l }\nonumber\\
&&\hspace*{1cm}=
[S(\tilde{x}^2_{\l})\ot S(\tilde{x}^1_{\l})\ot 1_{\mf{B}}]
[f\ot 1_{\mf{B}}](\Delta \ot id_{\mf{B}})(\tilde{q}_{\l}
\l (\tilde{x}^3_{\l})).\label{tql} 
\end{eqnarray}

Bicomodule algebras where introduced by Hausser and Nill 
in \cite{hn1}, under the name ``quasi-commuting pair of $H$-coactions". 

\begin{definition}\delabel{6.1}
Let $H$ be a quasi-bialgebra. An $H$-bicomodule algebra $\mb{A}$ 
is a quintuple $(\mb{A},\l, \r , \Phi_{\l }, \Phi_{\r }, \Phi _{\l , \r })$, 
where $\l $ and $\r $ are left and right $H$-coactions on $\mb {A}$, 
and where $\Phi _{\l}\in H\ot H\ot \mb {A}$, $\Phi _{\r }\in \mb {A}\ot H\ot H$ 
and $\Phi _{\l , \r }\in H\ot \mb {A}\ot H$ are invertible elements, such that
\begin{itemize}
\item[-] $(\mb {A}, \l , \Phi _{\l })$ is a left $H$-comodule algebra,
\item[-] $(\mb {A}, \r , \Phi _{\r })$ is a right $H$-comodule algebra,
\item[-] the following compatibility relations hold, 
for all $u\in \mb {A}$:
\begin{eqnarray}
&&\hspace*{-1cm}\Phi _{\l , \r }(\l \ot id)(\r (u))
=(id\ot \r )(\l (u))\Phi _{\l, \r }\label{bca1}\\
&&\hspace*{-1cm}(1_H\ot \Phi _{\l , \r })(id\ot \l \ot id)(\Phi _{\l , \r })
(\Phi _{\l }\ot 1_H)\nonumber\\
&&\hspace*{1cm}=(id\ot id\ot \r )(\Phi _{\l })
(\Delta \ot id\ot id)(\Phi _{\l , \r }) \label{bca2}\\
&&\hspace*{-1cm}(1_H\ot \Phi _{\r })(id\ot \r \ot id)(\Phi _{\l ,\r })
(\Phi _{\l , \r }\ot 1_H)\nonumber\\
&&\hspace*{1cm}= (id\ot id\ot \Delta )(\Phi _{\l , \r })
(\l \ot id\ot id) (\Phi _{\r }).\label{bca3}
\end{eqnarray}
\end{itemize}
\end{definition}

It was pointed out in \cite{hn1} that the following additional
relations hold in an $H$-bicomodule algebra $\mb {A}$:
$$
(id_H\ot id_{\mb {A}}\ot \va )(\Phi _{\l , \r })=1_H\ot 1_{\mb
{A}}, \mbox{${\;\;}$} (\va \ot id_{\mb {A}}\ot id_H)(\Phi _{\l ,
\r })= 1_{\mb {A}} \ot 1_H.
$$
As a first example, take $\mb {A}=H$, $\l =\r =\Delta $ and $\Phi
_{\l }=\Phi _{\r }= \Phi _{\l , \r }=\Phi $.\\
Let $(\mb{A},\l, \r , \Phi_{\l }, \Phi_{\r }, \Phi
_{\l , \r })$ be an $H$-bicomodule algebra; it is not hard to show that 
\begin{itemize}
\item[-] $(\mb{A}, \r \circ \tau , \l \circ \tau , (\Phi_{\r}^{-1})^{321}, 
(\Phi _{\l}^{-1})^{321}, (\Phi _{\l ,\r}^{-1})^{321})$ is an 
$H^{\rm cop}$-bicomodule algebra,  
\item[-] $(\mb{A}^{\rm op}, \r \circ \tau , \l \circ \tau , \Phi _{\r}^{321}, 
\Phi _{\l}^{321}, \Phi _{\l , \r}^{321})$ is an $H^{\rm op,cop}$-bicomodule 
algebra,
\item[-] $(\mb{A}^{\rm op}, \l , \r , \Phi_{\l}^{-1}, 
\Phi _{\r}^{-1}, \Phi _{\l ,\r}^{-1})$ is an $H^{\rm op}$-bicomodule algebra.
\end{itemize} 

We will use  the following
notation: 
$$
\Phi _{\l , \r }=\Theta ^1\ot \Theta ^2\ot \Theta ^3=
\ov {\Theta }^1\ot \ov {\Theta }^2\ot \ov {\Theta }^3;
$$
$$
\Phi ^{-1}_{\l , \r }=\theta ^1\ot \theta ^2\ot \theta ^3=
\ov {\theta }^1\ot \ov {\theta }^2\ot \ov {\theta }^3. 
$$

Let $H$ be a quasi-Hopf algebra and ${\mb A}$ an $H$-bicomodule algebra. 
We define two left $H\otimes H^{\rm op}$-coactions    
$\lambda _1, \lambda _2:{\mb A} \rightarrow (H\ot H^{\rm op})\ot {\mb A}$ on ${\mb A}$,
as follows:
\begin{eqnarray*}
\lambda _1 (u)&=&\left(u_{<0>_{[-1]}}\ot S^{-1}(u_{<1>})\right)\ot u_{<0>_{[0]}},\\
\lambda _2 (u)&=&\left(u_{[-1]}\ot S^{-1}(u_{[0]_{<1>}})\right)\ot u_{[0]_{<0>}},
\end{eqnarray*}
for all $u\in {\mb A}$. We also consider the following elements
$\Phi _{\lambda _1}, 
\Phi _{\lambda _2}\in (H\ot H^{\rm op})\ot 
(H\ot H^{\rm op})\ot {\mb A}$:
\begin{eqnarray*}
&&\Phi _{\l _1}=\left(\Theta ^1\tilde {X}^1_{\lambda }(\tilde {x}^1_{\rho })_{[-1]_1}
\ot S^{-1}(\tilde {x}^3_{\rho}g^2)\right)\\
&&\hspace*{2cm}
\ot \left(\Theta ^2_{[-1]}
\tilde {X}^2_{\lambda } (\tilde {x}^1_{\rho })_{[-1]_2}\ot 
S^{-1}(\Theta ^3\tilde {x}^2_{\rho }g^1)\right)\ot 
\Theta ^2_{[0]}\tilde {X}^3_{\lambda }(\tilde {x}^1_{\rho })_{[0]},\\
&&\Phi _{\l _2}=
\left(\tilde {Y}^1_{\lambda }\ot S^{-1}(\theta ^3\tilde {y}^3_{\rho}
(\tilde {Y}^3_{\lambda })_{<1>_2}g^2)\right)\\
&&\hspace*{2cm}
\ot \left(\theta ^1\tilde {Y}^2_{\lambda }\ot S^{-1}(\theta ^2_{<1>}
\tilde {y}^2_{\rho}(\tilde {Y}^3_{\lambda })_{<1>_1}g^1)\right)\ot 
\theta ^2_{<0>}\tilde {y}^1_{\rho }(\tilde {Y}^3_{\lambda })_{<0>}.
\end{eqnarray*}
It was proved in \cite{bpv3} that $(\mb{A}, \l _1, \Phi _{\l _1})$ and 
$(\mb{A}, \l _2, \Phi _2)$ are twist equivalent left $H\ot H^{\rm op}$-comodule 
algebras. In particular, if $H$ is a quasi-Hopf algebra then the notion 
of $H$-bicomodule algebra can be restated in terms of monoidal categories. 
In \seref{3} we will see that $\mb{A}$ can be also viewed in two twist 
equivalent ways as a right $H^{\rm op}\ot H$-comodule algebra. 

%%%%%%%%%%
\section{Doi-Hopf modules and rationality properties}\selabel{2}
%%%%%%%%%%
\subsection{Doi-Hopf modules}\selabel{2.1}
\setcounter{equation}{0}
The category of left (right) modules over a quasi-bialgebra is monoidal.
A coalgebra in ${}_H\Mm$ (resp. $\Mm_H$) is called a left (right)
$H$-module coalgebra. Thus a left $H$-module coalgebra is a left $H$-module
$C$ together with a comultiplication $\und :\ C\ra C\ot C$
and a counit $\une :\ C\ra k$ such that 
\begin{eqnarray}
&&\Phi (\und \ot id_C)(\und (c))=(id_C\ot \und )(\und (c)),\label{lmc1}\\
&&\und (h\cd c)=h_1\cd \una \ot h_2\cd \unb ,\label{lmc2}\\
&&\une (h\cd c)=\va (h)\une (c),
\end{eqnarray}
for all $c\in C$ and $h\in H$. Similarly, a right $H$-module 
coalgebra $C$ is a right $H$-module 
together with a comultiplication $\und :\ C\ra C\ot C$ and a counit $\une
:\ C\ra k$, satisfying the following relations
\begin{eqnarray}
&&(\und \ot id_C)(\und (c))\Phi ^{-1}=(id_C\ot \und )(\und (c)),\label{rmc1}\\
&&\und (c\cd h)=\una \cd h_1\ot \unb \cd h_2,\label{rmc2}\\
&&\une (c\cd h)=\une (c)\va (h),\label{rmc3}
\end{eqnarray}
for all $c\in C$ and $h\in H$. Here we used the Sweedler-type notation
$$
\und (c)=\una \ot
\unb , \mbox{${\;\;}$} (\und \ot id_C)(\und (c))=\unaa \ot
\unab \ot \unb \mbox{${\;\;}$ etc.}
$$
It is easy to see that  a left $H$-module coalgebra $C$ 
is in a natural way a right $H^{\rm op}$-module coalgebra (and vice versa). 

Let $H$ be a quasi-bialgebra and 
$C$ a right $H$-module coalgebra. A left $[C, H]$-Hopf module 
is a left $C$-comodule in the 
monoidal category ${\cal M}_H$. This definition was generalized 
in \cite{bc3}. 

\begin{definition}\delabel{2.1}
Let $H$ be a quasi-bialgebra over a field $k$, $C$ a right
$H$-module coalgebra and $(\mf {B}, \l , \Phi _{\l })$ a left
$H$-comodule algebra. A right-left $(H, \mf {B}, C)$-Hopf module
(or Doi-Hopf module) is a $k$-module $M$, with the following additional
structure: $M$ is right $\mf {B}$-module (the right action of
$\mfb $ on $m$ is denoted by $m\cd \mfb$), and we have a $k$-linear map
$\l _M:\ M\ra C\ot M$, such that the following relations hold, for
all $m\in M$ and $\mfb \in \mf {B}$:
\begin{eqnarray}
&&(\und \ot id_M)(\l _M(m))=(id_C\ot \l _M) (\l _M(m))\Phi _{\l },
\label{dhm1}\\
&&(\une \ot id_M)(\l _M(m))=m,\label{dhm2}\\
&&\l _M(m\cd \mfb )=m_{\{-1\}}\cd \mfb _{[-1]}\ot m_{\{0\}}\cd 
\mfb _{[0]}.\label{dhm3}
\end{eqnarray}
As usual, we use the Sweedler-type notation $\l _M(m)=m_{\{-1\}}\ot
m_{\{0\}}$. $^C{\cal M}(H)_{\mf {B}}$ is the category of right-left $(H,
{\mf B}, C)$-Hopf modules and right ${\mf B}$-linear, left
$C$-colinear $k$-linear maps.
\end{definition}

Let $M$ be a right $\mf{B}$-module; then  $C\ot M$ is
a right-left $(H, {\mf B}, C)$-Hopf module, with structure maps 
given by the following formulas
\begin{eqnarray}
&&(c\ot m)\cd \mfb =c\cd \mfb _{[-1]}\ot m\cd \mfb _{[0]},\label{sdhm1}\\
&&\l _{C\ot M}(c\ot m)=\una \cd \tilde{x}^1_{\l }\ot \unb \cd 
\tilde{x}^2_{\l}\ot m\cd \tilde{x}^3_{\l},\label{sdhm2}
\end{eqnarray} 
for all $c\in C$, $\mfb \in \mf{B}$ and $m\in M$. We obtain 
a functor ${\cal F}=C\ot \bullet :\  {\cal M}_{\mf{B}}\ra 
{}^C{\cal M}(H)_{\mf{B}}$.  The functor ${\cal F}$ sends a 
morphism $\vartheta $ to $id_C\ot \vartheta$. In particular,  
$C\ot \mf{B}\in {}^C{\cal M}(H)_{\mf{B}}$, via 
the structure maps
\begin{eqnarray*}
&&(c\ot \mfb)\cd \mfb ' =c\cd \mfb '_{[-1]}\ot \mfb \mfb '_{[0]},\\
&&\l _{C\ot \mf{B}}(c\ot \mfb )=\una \cd \tilde{x}^1_{\l }\ot \unb \cd 
\tilde{x}^2_{\l}\ot \mfb \tilde{x}^3_{\l},
\end{eqnarray*} 
for all $c\in C$ and $\mfb , \mfb '\in \mf{B}$. 

The functor 
${\cal F}$ has a left and a right adjoint, so it is an exact functor.

\begin{proposition}\prlabel{2.2}
Let $H$ be a quasi-bialgebra, $\mf{B}$ a left $H$-comodule algebra 
and $C$ a right $H$-module coalgebra. Then the functor 
${\cal F}=C\ot \bullet $ is a right adjoint of the forgetful functor 
$${}^C{\cal U}: \ {}^C{\cal M}(H)_{\mf{B}}\ra {\cal M}_{\mf{B}},$$
and a left 
adjoint of the functor 
$$\Hom_{\mf{B}}^C(C\ot \mf{B}, \bullet ): 
\ {}^C{\cal M}(H)_{\mf{B}}\ra {\cal M}_{\mf{B}}$$
defined as follows. For  $M\in {}^C{\cal M}(H)_{\mf{B}}$, 
$\Hom_{\mf{B}}^C(C\ot \mf{B}, M)$ is a right $\mf{B}$-module 
via the formula
$$
(\eta \cd \mfb )(c\ot \mfb ')=\eta (c\ot \mfb \mfb '),$$
for all $\eta \in \Hom_{\mf{B}}^C(C\ot \mf{B}, M)$, $c\in C$ and 
$\mfb , \mfb '\in \mf{B}$.
 For a morphism $\kappa :\ M\ra N$ in ${}^C{\cal M}(H)_{\mf{B}}$, we let
$$
\Hom_{\mf{B}}^C(C\ot \mf{B},\kappa )(\upsilon )=\kappa \circ \upsilon,$$
for all $\upsilon \in \Hom_{\mf{B}}^C(C\ot \mf{B}, M)$.   
\end{proposition}

\begin{proof}
Let $M$ be a right-left $(H, \mf{B}, C)$-Hopf module and $N$ a 
right $\mf{B}$-module. Define
$$\xi _{M, N}:\ \Hom_{\mf{B}}(M, N)\ra \Hom_{\mf{B}}^C(M, C\ot N),~~
\xi_{M, N}(\varsigma )(m)=m_{\{-1\}}\ot \varsigma (m_{\{0\}}),$$
for all $\varsigma \in \Hom_{\mf{B}}(M, N)$ and $m\in M$,
and
$$\zeta _{M, N}:\ \Hom_{\mf{B}}^C(M, C\ot N)\ra \Hom_{\mf{B}}(M, N),~~
\zeta _{M, N}(\chi )(m)=(\une \ot id_N)(\chi (m)),$$
for all $\chi \in \Hom_{\mf{B}}^C(M, C\ot N)$ and $m\in M$.
It is not hard to see that $\xi _{M, N}$ and $\zeta _{M, N}$ are 
well-defined natural transformations that are inverse to each other. \\
For $M\in {\cal M}_{\mf{B}}$ and 
$N\in {}^C{\cal M}(H)_{\mf{B}}$, we define 
$$\xi'_{M,N}:\ \Hom_{\mf{B}}^C(C\ot M, N)\to 
\Hom_{\mf{B}}(M, \Hom_{\mf{B}}^C(C\ot \mf{B}, N))$$
by
$$\xi '_{M, N}(\varsigma ')(m)(c\ot \mfb )=\varsigma '(c\ot m\cd \mfb )$$
and
$$\zeta '_{M, N}:\ \Hom_{\mf{B}}(M, \Hom_{\mf{B}}^C(C\ot \mf{B}, N))\to 
\Hom_{\mf{B}}^C(C\ot M, N)$$
by
$$\zeta '_{M, N}(\chi ')(c\ot m)=\chi '(m)(c\ot 1_{\mf{B}}),$$
for all $\varsigma '\in Hom_{\mf{B}}^C(C\ot M, N)$, 
$\chi '\in Hom_{\mf{B}}(M, Hom_{\mf{B}}^C(C\ot \mf{B}, N))$, $m\in M$, 
$c\in C$ and $\mfb \in \mf{B}$. Then $\xi'$ and $\zeta'$ are well-defined
natural transformations that are inverse to each other.
\end{proof}

\subsection{Doi-Hopf modules and comodules over a coring}\selabel{2.2}
It was proved in \cite{bc3} that the category of right-left Doi-Hopf modules is 
isomorphic to the category of right comodules over a suitable ${\mf{B}}$-coring.
For a general treatment of the theory of corings, we refer to 
\cite{Brzezinski02,BrzezinskiWisbauer}.

Let $R$ be a ring with unit.
An $R$-coring ${\cal C}$ is an $R$-bimodule 
together with two $R$-bimodule maps 
$\Delta _{\cal C}:\ {\cal C}\ra {\cal C}\ot _R{\cal C}$ 
and $\va _{\cal C}:\ {\cal C}\ra R$ such that the usual 
coassociativity and counit properties hold. A right ${\cal C}$-comodule 
is a right $R$-module together with a right $R$-linear map 
$\r ^r_M:\ M\ra M\ot _R{\cal C}$ such that 
\begin{eqnarray*}
&&(\r ^r_M\ot _Rid_{\cal C})\circ \r ^r_M=
(id_M\ot _R\Delta _{\cal C})\circ \r ^r_M,\\
&&(id_M\ot _R\va _{\cal C})\circ \r ^r_M=id_M.
\end{eqnarray*}  
A map $\hbar :\ M\ra N$ between two right ${\cal C}$-comodules is called 
right ${\cal C}$-colinear if $\hbar$ is right $R$-linear and 
$\r ^r_N\circ \hbar =(\hbar \ot _Rid_{\cal C})\circ \r ^r_M$. 
${\cal M}^{\cal C}$ will be the category of right ${\cal C}$-comodules 
and right ${\cal C}$-comodule maps. It is well-known that  ${\cal M}^{\cal C}$ 
is a Grothendieck category (in particular with enough injective objects) if ${\cal C}$ 
is flat as a left $R$-module. 

The category  ${}^{\cal C}{\cal M}$ of left ${\cal C}$-comodules 
and left ${\cal C}$-comodule maps can be introduced in a similar way.

Let $H$ be a quasi-bialgebra, $(\mf{B}, \l , \Phi _{\l })$ a left 
$H$-comodule algebra and $C$ a right $H$-module coalgebra. 
It was proved in \cite{bc3} that ${\cal C}=\mf{B}\ot C$ is a 
$\mf{B}$-coring. More precisely, ${\cal C}$ is a $\mf{B}$-bimodule via 
$$
\mfb \cd (\mfb '\ot c)=\mfb \mfb '\ot c~~{\rm and}~~
(\mfb \ot c)\cd \mfb '=\mfb \mfb '_{[0]}\ot c\cd \mfb '_{[-1]},
$$
and the comultiplication and counit are given by
$$
\Delta _{\cal C}(\mfb \ot c)=
(\mfb \tilde{x}^3_{\l}\ot \unb \cd \tilde{x}^2_{\l})\ot _{\mf{B}}
(1_{\mf{B}}\ot \una \cd \tilde{x}^1_{\l })$$
and
$$\va _{\cal C}(c\ot \mfb )=\une (c)\mfb ,$$
for all $\mfb , \mfb '\in \mf{B}$ and $c\in C$. Then we have an isomorphism of
categories (see \cite[Theorem 5.4]{bc3})
$${}^C{\cal M}(H)_{\mf{B}}\cong \Mm^\Cc.$$
${\cal C}$ is  free, and therefore flat, as a left $\mf{B}$-module, 
and we conclude that ${}^C{\cal M}(H)_{\mf{B}}\cong \Mm^\Cc$
is a Grothendieck category, and has therefore enough injectives
(see \cite[Prop. 1.2]{elkaoutit} or \cite[18.14]{BrzezinskiWisbauer}.

\subsection{Doi-Hopf modules and the generalized smash product}\selabel{2.3}
We now want to discuss when the category of Doi-Hopf modules is isomorphic
to a module category. In the case where $C$ is finite dimensional, this was
already done in \cite[Proposition 5.2]{bc3}. Let us explain this more precisely. 

Let $H$ be a quasi-bialgebra. A left $H$-module algebra $A$ is an 
algebra in the monoidal category  ${}_H{\cal M}$. This means that
$A$ is a left $H$-module with a multiplication $A\ot A\to A$ and a unit
element $1_A$  satisfying the 
following conditions: 
\begin{eqnarray}
&&(a a')a{''}=(X^1\cd a)[(X^2\cd a')(X^3\cd
a{''})],\label{ma1}\\
&&h\cd (a a')=(h_1\cd a)(h_2\cd a'),\label{ma2}\\
&&h\cd 1_A=\va (h)1_A,\label{ma3}
\end{eqnarray}
for all $a, a{'}, a{''}\in A$ and $h\in H$. 
Following  \cite{bc3}, we can define the generalized smash product of
a left $H$-module algebra $A$ and a left $H$-comodule algebra $\mf{B}$:
$A\gsm \mf{B}=A\ot B$ as a vector space, with multiplication
$$
(a\gsm \mfb )(a'\gsm \mfb ')=(\tilde{x}^1_{\l}\cd a)
(\tilde{x}^2_{\l}\mfb _{[-1]}\cd a')\gsm \tilde{x}^3_{\l}\mfb _{[0]}\mfb ',
$$    
for all $a, a'\in A$, $\mfb , \mfb '\in \mf{B}$. $A\gsm \mf{B}$ is an
associative algebra with unit $1_A\gsm 1_B$.

The linear dual $C^*$ of a right $H$-module coalgebra $C$
is a left $H$-module algebra. The multiplication is the convolution, the
unit is $\une$, and the left $H$-action is given by the formula
$(h\cd c^*)(c)=c^*(c\cd h)$, for all $h\in H$, $c^*\in C^*$ 
and $c\in C$. So we can consider the generalized smash product algebra 
$C^*\gsm \mf{B}$. Moreover, we have a functor 
$$\mf{G}: \ {}^C{\cal M}(H)_{\mf{B}}\ra {\cal M}_{C^*\gsm \mf{B}},
~~\mf{G}(M)=M,$$
with right $C^*\gsm \mf{B}$-action given by
$$
m\cd (c^*\gsm  \mfb )=c^*(m_{\{-1\}})m_{\{0\}}\cd \mfb ,
$$
for all $m\in M$, $c^*\in C^*$ and $\mfb \in \mf{B}$. If $C$ is finite 
dimensional, then 
$\mf{G}$ is an isomorphism of categories (see \cite{bc3}). 

Now let $C$ be infinite dimensional. We will show that 
the category ${}^C{\cal M}(H)_{\mf{B}}$ is 
%equivalent 
isomorphic 
to the category 
$\sigma _{C^*\gsm \mf{B}}[C\ot \mf{B}]$. Recall that if ${\cal A}$ is a 
Grothendieck category and $M$ is an object of ${\cal A}$ then 
$\sigma _{\cal A}[M]$ is the class of all objects  $N\in {\cal A}$ 
which are subgenerated by $M$, that is, $N$ is a subobject of a quotient
of directs sums of copies of
$M$. It is well-known that $\sigma _{\cal A}[M]$ is the 
smallest closed subcategory of ${\cal A}$ containing $M$, and that for any 
closed subcategory ${\cal D}$ of ${\cal A}$ there exists an object $M$ of 
${\cal A}$ such that ${\cal D}=\sigma _{\cal A}[M]$. In particular, 
$\sigma _{\cal A}[M]$ is a Grothendieck category, 
and has a fortiori enough injective objects. 
For more detail, the reader is invited to consult \cite{dnr, wis}.    

Let $H$ be a quasi-bialgebra, $C$ a right $H$-module coalgebra and $\mf{B}$ a 
left $H$-comodule coalgebra. For a right $C^*\gsm \mf{B}$-module $M$ we define the 
linear maps
$$\mu _M:\ M\ra \Hom(C^*\ot \mf{B}, M),~~
\mu _M(m)(c^*\ot \mfb )=m\cd (c^*\gsm \mfb ),$$
$$\nu _M:\ C\ot M\ra \Hom(C^*\ot \mf{B}, M),~~
\nu _M(c\ot m)(c^*\ot \mfb )=c^*(c)m\cd (\une \gsm \mfb ),$$
for all $m\in M$, $c\in C$, $c^*\in C^*$ and $\mfb \in \mf{B}$. It is easily
verified that $\mu _M$ and $\nu _M$ are injective linear maps. 

Inspired by \cite{dnt, doi} we propose the following 

\begin{definition}\delabel{2.3}
Let $M$ be a right $C^*\gsm \mf{B}$-module. We say that $M$ is rational if 
$\im(\mu _M)\subseteq \im(\nu _M)$.  $\Rat({\cal M}_{C^*\gsm \mf{B}})$ will be the 
full subcategory of ${\cal M}_{C^*\gsm \mf{B}}$ consisting of
rational $C^*\gsm \mf{B}$-modules. 
\end{definition} 

\begin{remark}
If $H$ is a coassociative bialgebra, $\mf{B}$ a left $H$-comodule algebra 
and $C$ a right $H$-module coalgebra 
in the usual sense, then a rational $C^*\gsm \mf{B}$-module $M$ is nothing else that 
a $C^*\gsm \mf{B}$-module wich is rational as a $C^*$-module. Here $M$ is 
viewed as a right $C^*$-module via the canonical algebra map 
$C^*\hookrightarrow C^*\gsm \mf{B}$. 
\end{remark}

It follows easily from \deref{2.3} that a right $C^*\gsm \mf{B}$-module is rational 
if and only if for every $m\in M$ there exist two finite sets $\{c_i\}_i\subseteq C$ 
and $\{m_i\}_i\subseteq M$ such that 
\begin{equation}\label{rat}
m\cd (c^*\gsm \mfb )=\sum \limits_i c^*(c_i)m_i\cd (\une \gsm \mfb ), 
\end{equation}
for all $c^*\in C^*$ and $\mfb \in \mf{B}$.
If $\{c'_j\}_j\subseteq C$ and $\{m'_j\}_j\subseteq M$ are two other finite 
sets satisfying (\ref{rat}), then 
$\sum \limits_i c_i\ot m_i=\sum \limits_j c'_j\ot m'_j$, because of 
the injectivity of the map $\nu _M$. So we have a well-defined map 
$$
\l _M:\ M\ra C\ot M,~~\l _M(m)=\sum \limits_i c_i\ot m_i,$$
for all $m\in M$.
If $C$ is finite dimensional, then any 
right $C^*\gsm \mf{B}$-module is rational. Indeed, we  take 
a finite dual basis
$\{(c_i, c^i)\}_i$ 
of $C$ and then consider for each $m\in M$ 
the finite sets $\{c_i\}_i\subseteq C$ and 
$\{m_i=m\cd (c^i\gsm 1_{\mf{B}})\}_i\subseteq M$. 

We now summarize the properties of rational $C^*\gsm \mf{B}$-modules. 

\begin{proposition}\prlabel{2.4}
Let $H$ be a quasi-bialgebra, $C$ a right $H$-module coalgebra and $\mf{B}$ a 
left $H$-comodule algebra. Then:
\begin{itemize}
\item[i)] A cyclic submodule of a rational $C^*\gsm \mf{B}$-module is finite 
dimensional. 
\item[ii)] If $M$ is a rational $C^*\gsm \mf{B}$-module and $N$ is a 
$C^*\gsm \mf{B}$-submodule of $M$, then $N$ and $M/N$ are rational 
$C^*\gsm \mf{B}$-modules.
\item[iii)] If $(M_i)_{i\in I}$ is a family of rational 
$C^*\gsm \mf{B}$-modules, then the direct sum $M=\bigoplus \limits _{i\in I}M_i$ 
in ${\cal M}_{C^*\gsm \mf{B}}$ is a rational $C^*\gsm \mf{B}$-module. 
\item[iv)] To any right $C^*\gsm \mf{B}$-module $M$ we can associate a unique 
maximal rational submodule $M^{\rm rat}$, which is equal to $\mu _M^{-1}(Im(\nu _M))$. 
It is also equal to the sum of all rational $C^*\gsm \mf{B}$-submodule of $M$. 
We have a left exact 
functor 
$$\Rat :\ {\cal M}_{C^*\gsm \mf{B}}\ra {\cal M}_{C^*\gsm \mf{B}},~~
\Rat(M)=   M^{\rm rat}.$$
\end{itemize}
\end{proposition}  

\begin{proof}
In fact this is a straightforward generalization of \cite[Theorem 2.1.3]{sw} and 
\cite[Theorem 2.2.6]{dnr}. Consider an element $m$ of 
a rational $C^*\gsm \mf{B}$-module $m$.
In what follows, $\{c_i\}_i\subseteq C$ and $\{m_i\}_i\subseteq M$ will then be two finite sets 
satisfying (\ref{rat}).

i) Let $m\cd (C^*\gsm \mf{B})$ be a cyclic submodule of a rational 
$C^*\gsm \mf{B}$-module $M$. $m\cd (C^*\gsm \mf{B})$ is generated by
the finite set $\{m_i\}_i$, so it is finite dimensional. 

ii) Take $m\in N\subset M$. Choose the $c_i$ in such a way that they are linearly
independent. 
Fix $j$, and take $c^*\in C^*$ such that $c^*(c_i)=\delta_{i,j}$. Then
$N$ contains 
$m\cd (c^*\gsm 1_{\mf{B}})=\sum \limits_i c^*(c_i)m_i\cd 
(\une \gsm 1_{\mf{B}})=m_j$, as needed.  

Let $\ov{m}$ be the class in $M/N$ represented by $m\in M$. For all
 $c^*\in C^*$ and $\mfb \in \mf{B}$, we have that $\ov{m}\cd (c^*\gsm \mfb )=
\sum \limits_i c^*(c_i)\ov{m_i}\cd (\une \gsm \mfb)$, 
and it follows that $M/N$ is a rational $C^*\gsm \mf{B}$-module.

iii) Every $m\in M$ can be written in a unique way as
$$m=\sum_{j\in J}m_j,$$
with $J\subset I$ finite, and $m_j\in M_j$.
Since $M_j$ is a rational $C^*\gsm \mf{B}$-module, there 
exist two finite sets $\{c^j_k\}_k\subseteq C$ and $\{m^j_k\}_k\subseteq M_j$ 
such that 
$$m_j\cd (c^*\gsm \mfb )=\sum \limits _kc^*(c^j_k)m^j_k\cd (\une \gsm \mfb ),$$
for all 
$c^*\in C^*$ and $\mfb \in \mf{B}$. We therefore have that
$$
m\cd (c^*\gsm \mfb )=\sum_{j\in J}m_j\cd (c^*\gsm \mfb )
=\sum \limits _{j\in J, k}c^*(c^j_k)m^j_k\cd (\une \gsm \mfb ),
$$
for all $c^*\in C^*$ and $\mfb \in \mf{B}$, and it follows that $M$ 
is a rational $C^*\gsm \mf{B}$-module.

iv) Let $M$ be a right $C^*\gsm \mf{B}$-module. We define 
$M^{\rm rat}=\mu ^{-1}_M(\im(\nu _M))$. \\
We first prove that $M^{\rm rat}$ is a right $C^*\gsm \mf{B}$-module. 
Take $m\in M^{\rm rat}$. Then there exist two finite sets $\{c_i\}_i\subseteq C$ 
and $\{m_i\}_i\subseteq M$ such that 
$m\cd (c^*\gsm \mfb )=c^*(c_i)m_i\cd (\une \gsm \mfb )$, 
for all $c^*\in C^*$ and $\mfb \in \mf{B}$. Therefore:
\begin{eqnarray*}
&&\hspace*{-1.5cm}
(m\cd (c^*\gsm \mfb ))\cd (d^*\gsm \mfb ')
=m\cd ((c^*\gsm \mfb )(d^*\gsm \mfb '))\\
&=&
m\cd ((\tilde{x}^1_{\l }\cd c^*)(\tilde{x}^2_{\l }\mfb _{[-1]}\cd d^*)\gsm \tilde{x}^3_{\l}
\mfb _{[0]}\mfb ')\\
&=&c^*((c_i)_{\un{1}}\cd \tilde{x}^1_{\l})
d^*((c_i)_{\un{2}}\cd \tilde{x}^2_{\l}\mfb _{[-1]})(m_i\cd 
(\une \gsm \tilde{x}^3_{\l}\mfb _{[0]}))\cd (\une \gsm \mfb '),   
\end{eqnarray*}
for all $m\in M$, $c^*, d^*\in C^*$ and $\mfb , \mfb '\in \mf{B}$. 
Thus $m\cd (c^*\gsm \mfb )\in M^{\rm rat}$, hence
$M^{\rm rat}$ is a $C^*\gsm \mf{B}$-submodule of $M$.  Using an argument 
similar to the one in the first part of the proof of assertion ii), we can easily check that 
$M^{\rm rat}$ is a rational $C^*\gsm \mf{B}$-module.

Let $N$ be a rational $C^*\gsm \mf{B}$-submodule of $M$, this means 
$\im(\mu _N)\subseteq \im(\nu _N)$. Then 
$$\mu _M(N)=\mu _N(N)\subseteq 
\im(\nu _N)\subseteq \im(\nu _M),$$
hence $N\subseteq \mu ^{-1}_M(\im(\nu _M))=M^{\rm rat}$, and
we conclude that $M^{\rm rat}$ is the unique maximal rational submodule of $M$. Assertions 
ii) and iii) show that $M^{\rm rat}$ is also 
equal to the sum of all rational $C^*\gsm \mf{B}$-submodules of $M$. 
The proof of the final assertion is identical 
to the proof of \cite[Theorem 2.2.6 iv)]{dnr}.     
\end{proof}

We are now able to prove the main result of this Section.

\begin{theorem}\thlabel{2.5}
Let $H$ be a quasi-bialgebra, $\mf{B}$ a left $H$-comodule algebra 
and $C$ a right $H$-module coalgebra. 
The categories ${}^C{\cal M}(H)_{\mf{B}}$ and 
$\Rat({\cal M}_{C^*\gsm \mf{B}})$ are isomorphic, and 
$\Rat({\cal M}_{C^*\gsm \mf{B}})$ is equal to 
$\sigma _{C^*\gsm \mf{B}}[C\ot \mf{B}]$.  
\end{theorem}

%\begin{remark}
%It is well-known that the category $\sigma[M]$ is a Grothendieck category,
%so it follows from \thref{2.5} that
%${}^C{\cal M}(H)_{\mf{B}}$ is a Grothendieck category, an observation that
%we already made before.
%\end{remark}

\begin{proof}
Recall that we have a functor
$$\mf{G}: {}^C{\cal M}(H)_{\mf{B}}\ra {\cal M}_{C^*\gsm \mf{B}},~~
\mf{G}(M)=M,$$
with right $C^*\gsm \mf{B}$-action given by the formula
$$m\cd (c^*\gsm \mfb )=c^*(m_{\{-1\}})m_{\{0\}}\cd \mfb .$$
It is clear that $\mf{G}(M)$ is rational as a $C^*\gsm \mf{B}$-module.

Consider a rational $C^*\gsm \mf{B}$-module $M$. Then
$M$ is a right $\mf{B}$-module through $m\cd \mfb =m\cd (\une \gsm \mfb )$.
We define a linear map $\l _M:\ M\ra C\ot M$ as follows: 
$$
\l _M(m)=\sum_i c_i\ot m_i$$
if and only if
$$m\cd (c^*\gsm \mfb )=\sum \limits_i c^*(c_i)m_i\cd (\une \gsm \mfb ),$$
for all
$c^*\in C^*$ and $\mfb \in \mf{B}$. It is clear that $\l _M$ is well-defined.

Fix $i$, and assume that
$$\l _M(m_i)=\sum_j c^i_j\ot m^i_j,$$
or, equivalently,
$$m_i\cd (c^*\gsm \mfb )=c^*(c^i_j)m^i_j\cd (\une \gsm \mfb),$$
for all $c^*\in C^*$ and $\mfb \in \mf{B}$. Therefore
we have, for all $c^*, d^*\in C^*$ and $m\in M$, that
\begin{eqnarray*}
&&\hspace*{-1cm}
\le c^*\ot d^*\ot id_M, (id_C\ot \l _M)(\l _M(m))\cd \Phi _{\l}\ri 
=c^*(c_i\cd \tilde{X}^1_{\l})d^*(c^i_j\cd 
\tilde{X}^2_{\l})m^i_j\cd \tilde{X}^3_{\l}\\
&&=c^*(c_i\cd \tilde{X}^1_{\l})m_i\cd (\tilde{X}^2\cd d^*\gsm \tilde{X}^3_{\l})
=(m\cd (\tilde{X}^1_{\l}\cd c^*\gsm 1_{\mf{B}}))\cd 
(\tilde{X}^2_{\l}\cd d^*\gsm \tilde{X}^3_{\l})\\
&&=m\cd ((\tilde{X}^1_{\l}\cd c^*\gsm 1_{\mf{B}})(\tilde{X}^2_{\l}\cd 
d^*\gsm \tilde{X}^3_{\l}))
=m\cd (c^*d^*\gsm 1_{\mf{B}})\\
&&=c^*d^*(c_i)m_i=\le c^*\ot d^*\ot id_M, (\und \ot id_M)(\l _M(m))\ri ,   
\end{eqnarray*} 
proving that (\ref{dhm1}) is satisfied.  (\ref{dhm2}) is trivial since 
$m=m\cd (\une \gsm 1_{\mf{B}})$ for all $m\in M$. We will next
prove that (\ref{dhm3}) holds. First observe that
\begin{eqnarray*}
&&\hspace*{-2cm}
(m\cd \mfb )\cd (c^*\gsm \mfb ')=(m\cd (\une \gsm \mfb ))\cd (c^*\gsm \mfb ')
=m\cd ((\une \gsm \mfb )(c^*\gsm \mfb '))\\
&=&m\cd (\mfb _{[-1]}\cd c^*\gsm \mfb _{[0]}\mfb ')=c^*(c_i\cd \mfb _{[-1]})
(m_i\cd \mfb _{[0]})\cd (\une \gsm \mfb '), 
\end{eqnarray*}
for all $m\in M$, $c^*\in C^*$ and $\mfb , \mfb '\in \mf{B}$.
This shows 
that $\l _M(m\cd \mfb )=c_i\cd \mfb _{[-1]}\ot m_i\cd \mfb _{[0]}$, as needed,
and it follows that $M\in {}^C{\cal M}(H)_{\mf{B}}$.

Let $\eta :\ M\ra N$ be a morphism of rational $C^*\gsm \mf{B}$-modules.
Take $m\in M$, and assume that $\l _M(m)=\sum \limits_i c_i\ot m_i$. We compute
that
$$
\eta (m)\cd (c^*\gsm \mfb )=\eta (m\cd (c^*\gsm \mfb ))
=\sum \limits_ic^*(c_i)\eta (m_i)\cd (\une \gsm \mfb ),
$$
for all $c^*\in C^*$ and $\mfb \in \mf{B}$. This is equivalent to
$\l _N(\eta (m))=c_i\ot \eta (m_i)=(id_C\ot \eta )(\l _M(m))$, 
hence $\eta $ is left $C$-colinear. It is clear that $\eta $ is right $\mf{B}$-linear,
so $\eta$ is a morphism in  $ {}^C{\cal M}(H)_{\mf{B}}$, and we have a functor
$\mf{F}: \Rat({\cal M}_{C^*\gsm \mf{B}})\ra {}^C{\cal M}(H)_{\mf{B}}$, 
which is inverse to $\mf{G}$.\\

The proof of the fact that 
%${}^C{\cal M}(H)_{\mf{B}}$ 
$\Rat({\cal M}_{C^*\gsm \mf{B}})$
and 
$\sigma _{C^*\gsm \mf{B}}[C\ot \mf{B}]$ are 
%equivalent categories 
equal is similar to 
the proof of \cite[Lemma 3.9]{mz}. 
%It suffices to show that the image 
%of the functor $\mf{G}$ is $\sigma_{C^*\gsm \mf{B}}[C\ot \mf{B}]$.

We first show that 
%$\mf{G}(M)\in \sigma _{C^*\gsm \mf{B}}[C\ot \mf{B}]$,
$M\in \sigma _{C^*\gsm \mf{B}}[C\ot \mf{B}]$ 
for every 
%$M\in  {}^C{\cal M}(H)_{\mf{B}}$. 
$M\in \Rat({\cal M}_{C^*\gsm \mf{B}})$. 
We recall first that a 
right $C^*\gsm \mf{B}$-module belongs to $\sigma _{C^*\gsm \mf{B}}[C\ot \mf{B}]$ 
if and only if there exists a set $I$, a right $C^*\gsm \mf{B}$-module $N$, and 
two $C^*\gsm \mf{B}$-linear maps $\iota : M\ra N$ and 
$\pi : (C\ot \mf{B})^{(I)}\ra N$ such that $\iota $ is injective and $\pi $ 
is surjective. 

Let $M$ be a rational $C^*\gsm \mf{B}$-module. Thus 
$M\in {}^C{\cal M}(H)_{\mf{B}}$ and $\mf{G}(M)=M$ as right $C^*\gsm \mf{B}$-modules. 
It is easy to check that $\iota =\l _M:\ M\ra C\ot M$ 
is an injective morphism in ${}^C{\cal M}(H)_{\mf{B}}$; here 
$C\ot M$ has the right-left Doi-Hopf module structure given 
by (\ref{sdhm1}, \ref{sdhm2}). The map 
$$
\pi _M:\ \mf{B}^{(M)}\ra M,~~\pi _M((\mfb _m)_m)=m\cd \mfb _m,
$$   
is a surjective right $\mf{B}$-linear map. By \prref{2.2}, it provides 
a surjective morphism $id_C\ot \pi _M:\ C\ot \mf{B}^{(M)}\ra C\ot M$ in 
${}^C{\cal M}(H)_{\mf{B}}$, and since $C\ot \mf{B}^{(M)}\cong (C\ot \mf{B})^{(M)}$ 
in ${}^C{\cal M}(H)_{\mf{B}}$, we conclude that there exists a surjective morphism 
$\pi :\ (C\ot \mf{B})^{(M)}\ra C\ot M$ in ${}^C{\cal M}(H)_{\mf{B}}$, so 
$M\in \sigma _{C^*\gsm \mf{B}}[C\ot \mf{B}]$. 

Take $M\in \sigma _{C^*\gsm \mf{B}}[C\ot \mf{B}]$. Then we have
 right $C^*\gsm \mf{B}$-morphisms
$\iota :\ M\ra N$ and $\pi :\ (C\ot \mf{B})^{(I)}\ra N$ 
such that $\iota$ is injective and $\pi$ is surjective. 
Since the right $C^*\gsm \mf{B}$-module $C\ot \mf{B}$ lies
in the image of $\mf{G}$, it follows that $C\ot \mf{B}$ is a rational module. 
By \prref{2.4} iii), $(C\ot \mf{B})^{(I)}$ is a rational module too and since 
$\pi $ is surjective we deduce from \prref{2.4} ii) that $N$ is rational. Finally, 
from \prref{2.4} ii) and the fact that
$\iota $ is injective, it follows that $M$ is a rational 
$C^*\gsm \mf{B}$-module. 
%From the first part of the proof, it follows that $M$ is a right-left 
%Doi-Hopf module, and that $M=\mf{G}(M)$. 
\end{proof}

\begin{remark}
It is well-known that the category $\sigma _{\cal A}[M]$ is a Grothendieck category,
so it follows from \thref{2.5} that
${}^C{\cal M}(H)_{\mf{B}}$ is a Grothendieck category, an observation that
we already made before.
\end{remark}

\begin{corollary}\cite[Proposition 5.2]{bc3}\colabel{2.6}
Let $H$ be a quasi-bialgebra, $\mf{B}$ a left $H$-comodule algebra 
and $C$ a finite dimensional right $H$-module coalgebra. Then the categories 
${}^C{\cal M}(H)_{\mf{B}}$ and ${\cal M}_{C^*\gsm \mf{B}}$ are isomorphic. 
\end{corollary}

\begin{corollary}
Let $M$ be a right-left $(H, \mf{B}, C)$-Hopf module. Then the following 
assertions hold:
\begin{itemize}
\item[i)] The right-left $(H, \mf{B}, C)$-Hopf submodule generated by an element 
of $M$ is finite dimensional. 
\item[ii)] $M$ is the sum of its finite dimensional $(H, \mf{B}, C)$-Hopf 
submodules.  
\end{itemize}
\end{corollary}

\begin{proof}
i) $M$ is a rational $C^*\gsm \mf{B}$-module, so the right-left 
$(H, \mf{B}, C)$-Hopf submodule generated by an element $m\in M$ 
coincides with the cyclic $C^*\gsm \mf{B}$-submodule generated by $m$, 
by \thref{2.5}. We know from \prref{2.4} i) that it is finite dimensional.  

ii) View $M$ as a rational $C^*\gsm \mf{B}$-module. Obviously, 
$M$ is the sum of its cyclic $C^*\gsm \mf{B}$-submodules and all of these 
are finite dimensional right-left $(H, \mf{B}, C)$-Hopf submodules. 
\end{proof}

\subsection{Doi-Hopf modules and Koppinen's smash product}\selabel{2.4}
We begin this Section with some general results about corings, taken from
\cite[Sec. 19 and 20]{BrzezinskiWisbauer}. Let $R$ be a ring, and $\Cc$
an $R$-coring. Then ${}^*\Cc={}_R\Hom(\Cc,R)$ is an $R$-ring, with
multiplication
$$(\varphi\# \psi)(c)=\psi(c_{(1)}\varphi(c_{(2)})),$$
for all $\varphi,\psi\in {}^*\Cc$. We have a functor $F:\ \Mm^\Cc\to
\Mm_{{}^*\Cc}$, $F(M)=M$, with $m\cdot\varphi=m_{[0]}\varphi(m_{[1]})$,
for all $m\in M$ and $\varphi\in {}^*\Cc$. If $\Cc$ is finitely generated
and projective as a left $R$-module, then $F$ is an isomorphism of categories.\\
Assume now that $\Cc$ is locally projective as a left $R$-module.
Then $\Mm^\Cc$ is isomorphic to the category $\sigma_{{}^*\Cc}[\Cc]$,
see \cite[19.3]{BrzezinskiWisbauer}; observe that our multiplication on
${}^*\Cc$ is opposite to the one from \cite{BrzezinskiWisbauer}, so that
left ${}^*\Cc$-modules in \cite{BrzezinskiWisbauer} are our right
${}^*\Cc$-modules.\\
Take a right ${}^*\Cc$-module $M$. $m\in M$ is called rational if there
exists $\sum \limits_i m_i\ot c_i\in M\ot_R \Cc$ such that $m\cdot \varphi=
\sum \limits_i m_i\varphi(c_i)$, for all $\varphi\in {}^*\Cc$. Then $\Rat^\Cc(M)=
\{m\in M~|~m~{\rm~is~rational}\}$ is a right $\Cc$-comodule, and we obtain
a functor $\Rat^\Cc:\ \Mm_{{}^*\Cc}\to \Mm^\Cc$, which is right adjoint to $F$.
$M$ is called rational if $\Rat^\Cc(M)=M$.
$\Mm^\Cc$ is isomorphic to the full category of $\Mm_{{}^*\Cc}$ consisting
of rational right ${}^*\Cc$-modules.\\

Now let $H$ be a quasi-bialgebra, $\mf{B}$ a left $H$-comodule algebra,
and $C$ a right $H$-comodule algebra. We consider the $\mf{B}$-coring
$\Cc=\mf{B}\ot C$ from \seref{2.2}. Since we work over a field $k$, 
$C$ is projective as 
a $k$-module, hence $\Cc$ is projective (and a fortiori locally projective) 
as a left $\mf{B}$-module. Hence we can apply the above results to this 
situation. We have an isomorphism of vector spaces 
$${}^*\Cc={}_{\mf{B}}\Hom (\mf{B}\ot C, \mf{B})\cong \Hom (C, \mf{B}).$$
The multiplication on ${}^*\Cc$ can be transported to a multiplication on
$\Hom (C, \mf{B})$. This multiplication makes $\Hom (C, \mf{B})$ into a
$B$-ring $\# (C, \mf{B})$, which we will call the Koppinen smash product.
The multiplication is given by the following formula:
\begin{equation}\eqlabel{kop}
(f\# g)(c)=
\tilde{x}^3_{\lambda}f(c_{\ul{1}}\cdot \tilde{x}^1_{\lambda})_{[0]}
g\Bigl(
c_{\ul{2}}\cdot \tilde{x}^2_{\lambda}f(c_{\ul{1}}\cdot \tilde{x}^1_{\lambda})_{[1]}
\Bigr).
\end{equation}
In the situation where $H$ is an associative bialgebra, we recover the
smash product introduced first by Koppinen in \cite{koppinen}. The relation to
the generalized smash product introduced in \seref{2.3} is discussed in
\prref{2.4.1}.

\begin{proposition}\prlabel{2.4.1}
The $k$-linear map
$$\alpha:\ C^*\gsm \mf{B}\to \#(C,B),~~\alpha(c^*\gsm \mf{b})=f,$$
with $f(c)=\lan c^*,c\ran b$, for all $c\in C$, is a morphism of $\mf{B}$-rings.
It is an isomorphism if $C$ is finite dimensional.
\end{proposition}

\begin{proof}
We have to show that $\alpha$ is multiplicative. Take
$c^*\gsm \mf{b}, d^*\gsm \mf{b}'\in C^*\gsm \mf{B}$, and write
 $\alpha(c^*\gsm \mf{b})=f$, $\alpha(d^*\gsm \mf{b}')=g$.
Using \equref{kop}, we compute that
$$(f\# g)(c)=
\tilde{x}^3_{\lambda}\lan c^*, c_{\ul{1}}\cdot \tilde{x}^1_{\lambda}\ran \mf{b}_{[0]}
\lan d^*,c_{\ul{2}}\cdot \tilde{x}^2_{\lambda} \mf{b}_{[1]}\ran \mf{b}'.$$
We also have that
$$(c^*\gsm \mf{b})(d^*\gsm \mf{b}')=
(\tilde{x}^1_{\lambda}\cdot c^*)(\tilde{x}^2_{\lambda}\mf{b}_{[-1]}\cdot d^*)\gsm
\tilde{x}^3_{\lambda}\mf{b}_{[0]}\mf{b}',$$
so
$$
\alpha\bigl((c^*\gsm \mf{b})(d^*\gsm \mf{b}')\bigr)(c)
=\lan c^*, c_{\ul{1}}\cdot \tilde{x}^1_{\lambda}\ran
\lan d^*,c_{\ul{2}}\cdot \tilde{x}^2_{\lambda} \mf{b}_{[1]}\ran
\tilde{x}^3_{\lambda} \mf{b}_{[0]}\mf{b}',$$
as needed.
\end{proof}

Take a right $\#(C,B)$-module $M$. $m\in M$ is rational if there exists
$\sum \limits_i m_i\ot c_i\in M\ot C$ such that
$$m\cdot f=\sum \limits_i m_if(c_i),$$
for all $f:\ C\to B$. $M$ is called rational if every $m\in M$ is rational.

\begin{corollary}\colabel{2.4.2}
Now let $H$ be a quasi-bialgebra, $\mf{B}$ a left $H$-comodule algebra,
and $C$ a right $H$-comodule algebra. Then the category ${}^C\Mm(H)_{\mf{B}}$
is isomorphic to the full subcategory of $\Mm_{\#(C,B)}$, which is also
equal to 
$\sigma_{\#(C,\mf{B})}(\mf{B}\ot C)$.
\end{corollary}

\subsection{Left, right and right-left Doi-Hopf modules}\selabel{2.5}
For the sake of completeness, we also define the other Doi-Hopf module categories. 
In fact, we have four different types of Doi-Hopf modules. The first one was 
already studied, namely the right-left version. We also have the left-right,
right-right and left-left versions.

\begin{definition}\delabel{2.7}
Let $H$ be a quasi-bialgebra, $\mfA$ a right $H$-comodule algebra and 
$\mf{B}$ a left $H$-comodule algebra. In the statements below we assume  
in i) and iii) that $C$ is a left $H$-module coalgebra, and in ii) that 
$C$ is a right $H$-module coalgebra, respectively.\\

i) A left-right 
$(H, \mfA , C)$-Hopf module (or Doi-Hopf module) 
is a left $\mfA$-module $M$ together with 
a $k$-linear map $\r _M:\ M\ra M\ot C$, $\r _M(m)=m_{(0)}\ot m_{(1)}$,
such that the following relations hold, for all 
$m\in M$ and $\mfa \in \mfA$:
\begin{eqnarray}
&&\Phi _{\r }\cd (\r _M\ot id_C)(\r _M(m))=(id_M\ot \und )(\r _M(m)),\label{lrdhm1}\\
&&(id_M\ot \une )(\r _M(m))=m,\label{lrdhm2}\\
&&\r _M(\mfa \cd m)=\mfa _{\le0\ri}\cd m_{(0)}\ot \mfa _{\le1\ri}\cd m_{(1)}.\label{lrdhm3}
\end{eqnarray}
${}_{\mfA}{\cal M}(H)^C$ is the category of left-right $(H, \mfA , C)$-Hopf 
modules and left $\mfA$-linear, right $C$-colinear maps.\\

ii) A right-right 
$(H, \mfA , C)$-Hopf module (or Doi-Hopf module) is a right $\mfA$-module $M$ together with 
a $k$-linear map $\r _M:\ M\ra M\ot C$, $\r _M(m)=m_{(0)}\ot m_{(1)}$,
such that the following relations hold, for all 
$m\in M$ and $\mfa \in \mfA$:
\begin{eqnarray}
&&(\r _M\ot id_M)(\r _M(m))=(id_M\ot \und )(\r _M(m))\cd \Phi _{\r},\label{rrdhm1}\\
&&(id_M\ot \une )(\r _M(m))=m,\label{rrdhm2}\\
&&\r _M(m\cd \mfa )=m_{(0)}\cd \mfa _{\le0\ri}\ot m_{(1)}\cd \mfa _{\le1\ri}.\label{rrdhm3}
\end{eqnarray} 
${\cal M}(H)^C_{\mfA}$ is the category of right-right $(H, \mfA , C)$-Hopf 
modules and left $\mfA$-linear, right $C$-colinear maps.\\

iii) A left-left 
$(H, \mf{B}, C)$-Hopf module (or Doi-Hopf module) is a left $\mf{B}$-module $M$ 
together with a left $k$-linear map $\l _M:\ M\ra C\ot M$, 
$\l _M(m)=m_{\{-1\}}\ot m_{\{0\}}$, such that the following relations 
hold, for all $m\in M$ and $\mfb \in \mf{B}$: 
\begin{eqnarray}
&&\Phi _{\l}\cd (\und \ot id_M)(\l _M(m))=(id_C\ot \l _M)(\l _M(m)),\label{lldhm1}\\
&&(\une \ot id_M)(\l _M(m))=m,\label{lldhm2}\\
&&\l _M(\mfb \cd m)=\mfb _{[-1]}\cd m_{\{-1\}}\ot \mfb _{[0]}\cd m_{\{0\}}.\label{lldhm3}
\end{eqnarray}
${}_{\mf{B}}^C{\cal M}(H)$ is the category of left-left $(H, \mfA , C)$-Hopf 
modules and left $\mfA$-linear,left $C$-colinear maps.
\end{definition}

We remind that if $(\mf{A}, \r , \Phi _{\r})$ is a right $H$-comodule algebra 
then $\un{\mfA}^{\rm op}=(\mfA ^{\rm op}, \r \circ \tau _{\mf{A}, H}, (\Phi _{\r })^{321})$ is a 
left $H^{\rm op,cop}$-comodule algebra and 
$\un{\mfA}=(\mfA , \r \circ \tau _{\mf{A}, H}, (\Phi _{\r }^{-1})^{321})$ is a left 
$H^{\rm cop}$-comodule algebra, where $\tau _{\mf{A}, H}:\ \mfA \ot H\ra H\ot \mfA $ is the 
switch map. Also, if $\mf{B}$ is a left $H$-comodule algebra then 
$\mf{B}^{\rm op}=(\mf{B}^{\rm op}, \l , \Phi _{\l }^{-1})$ 
is a left $H^{\rm op}$-comodule algebra. 

On the other hand, if $C$ is a left $H$-module coalgebra then 
$C$ is a right $H^{\rm op}$-module coalgebra and $C^{\rm cop}$ is a 
right $H^{\rm op,cop}$-module coalgebra (and vice versa). So if $C$ is a right 
$H$-module coalgebra then $C^{\rm cop}$ is a right $H^{\rm cop}$-module coalgebra. 

Having these correspondences one can easily see that 
$$
{}_{\mfA}{\cal M}(H)^C\cong {}^{C^{\rm cop}}{\cal M}(H^{\rm op,cop})_{\un{\mfA}^{\rm op}},
~~{\cal M}(H)_{\mfA}^C\cong {}^{C^{\rm cop}}{\cal M}(H^{\rm cop})_{\un{\mfA}},
$$
$${}^C_{\mf{B}}{\cal M}(H)\cong {}^C{\cal M}(H^{\rm op})_{\mf{B}^{\rm op}}.$$
It follows that the four different types of Doi-Hopf modules are isomorphic to
categories of comodules over suitable corings, and they are 
Grothendieck categories with 
enough injective objects. On the other hand, if $C$ is finite dimensional then the 
above categories are isomorphic to categories of modules over certain generalized 
smash product algebras. More precisely, we have:

\begin{remarks}\reslabel{2.8}
i) The category ${}_{\mfA}{\cal M}(H)^C$ is isomorphic to the category 
of right modules over the generalized smash product 
$(C^{\rm cop})^*\gsm \un{\mfA}^{\rm op}$ (over $H^{\rm op,cop}$), 
and therefore also to the category 
of left modules over $((C^{\rm cop})^*\gsm \un{\mfA}^{\rm op})^{\rm op}$. Is not hard 
to see that the multiplication rule in $((C^{\rm cop})^*\gsm \un{\mfA}^{\rm op})^{\rm op}$ is 
$$
(c^*\gsm \mfa )(d^*\gsm \mfa ')=(c^*\cd \mfa '_{\le1\ri}\tilde{x}^2_{\r})
(d^*\cd \tilde{x}^3_{\r})\gsm \mfa \mfa '_{\le0\ri}\tilde{x}^1_{\r},
$$ 
for all $c^*, d^*\in C^*$ and $\mfa , \mfa '\in \mfA$, where 
$(c^*\cd h)(c)=c^*(h\cd c)$ for all $c^*\in C^*$, $h\in H$ and $c\in C$. 
Therefore, under the trivial permutation of tensor factors we have that 
$((C^{\rm cop})^*\gsm \un{\mfA}^{\rm op})^{\rm op}=\mfA \gtl C^*$, the right generalized 
smash product between the right $H$-comodule algebra $\mfA$ and the right 
$H$-module algebra $C^*$ (see \cite{bpv3} for more details). We conclude that
${}_{\mfA}{\cal M}(H)^C\cong 
{}_{\mfA \gtl C^*}{\cal M}$ if $C$ is finite dimensional.\\

ii) The above arguments entail that
$$
{\cal M}(H)_{\mfA}^C\cong {}^{C^{\rm cop}}{\cal M}(H^{\rm op})_{\un{\mfA}}\cong 
{\cal M}_{(C^{\rm cop})^*\gsm \un{\mfA}}\cong 
{}_{((C^{\rm cop})^*\gsm \un{\mfA})^{\rm op}}{\cal M},
$$
where the generalized smash product is over $H^{\rm cop}$. The explicit 
formula for the multiplication $\odot$ on $((C^{\rm cop})^*\gsm \un{\mfA})^{\rm op}$ is 
given by 
\begin{equation}\label{stgsm}
(c^*\gsm \mfa )\odot (d^*\gsm \mfa ')=(\tilde{X}^2_{\r}\mfa '_{\le1\ri}\cd c^*)
(\tilde{X}^3_{\r }\cd d^*)\gsm \tilde{X}^1_{\r }\mfa '_{\le0\ri}\mfa ,
\end{equation}
for all $c^*, d^*\in C^*$ and $\mfa , \mfa '\in \mf{A}$. 

iii) Obviously, ${}^C_{\mf{B}}{\cal M}(H)\cong {\cal M}_{C^*\gsm \mf{B}^{\rm op}}$, 
where the generalized smash product is taken over $H^{\rm op}$. 
\end{remarks}

Let $H$ be a quasi-Hopf algebra and $C$ a finite dimensional right 
$H$-module coalgebra. It was proved in \cite[Proposition 3.2]{bn} that 
the category ${\cal M}^C_H={\cal M}(H)^C_H$ is isomorphic to the category 
of left modules over the smash product algebra $C^*\# H=C^*\gsm H$. Now, 
by \resref{2.8} ii) the category ${\cal M}^C_H$ is isomorphic to the category 
of left modules over $((C^{\rm cop})^*\# \un{H})^{\rm op}$. The next result 
shows that the smash product algebras 
$C^*\# H$ and $((C^{\rm cop})^*\# \un{H})^{\rm op}$ are isomorphic.

\begin{proposition}\prlabel{2.9}
Let $H$ be a quasi-Hopf algebra and $C$ a right $H$-module coalgebra. Then 
the map 
$$\varphi :\ C^*\# H\ra ((C^{\rm cop})^*\# \un{H})^{\rm op},~~
\varphi (c^*\# h)=\smi (\mf{q}^1h_1g^1)\cd c^*\# \smi (\mf{q}^2h_2g^2)
$$
is an algebra isomorphism. Here 
$f^{-1}=g^1\ot g^2$ is the element defined by (\ref{g}) and 
$\tilde{q}_{\Delta }=\mf{q}^1\ot \mf{q}^2$ is the 
element $\tilde{q}_{\l}$ defined in (\ref{tpql}), in the special case where $\mf{B}=H$. 
\end{proposition}

\begin{proof}
We first show that $\varphi $ is an algebra map. Let $f=f^1\ot f^2$, 
$f^{-1}=g^1\ot g^2=G^1\ot G^2$ and 
$\tilde{q}_{\Delta }=\mf{q}^1\ot \mf{q}^2=\mf{Q}^1\ot \mf{Q}^2$ be the 
elements defined by (\ref{f}), (\ref{g}) and (\ref{tpql}), respectively. 
We compute 
\begin{eqnarray*}
&&\hspace*{-2cm}
\v ((c^*\# h)(d^*\# h'))=
\v ((x^1\cd c^*)(x^2h_1\cd d^*)\# x^3h_2h')\\[1mm]
&{{\rm (\ref{ma2}, \ref{ca})}\atop =}&
[\smi (f^2\mf{q}^1_2x^3_{(1, 2)}(h_2h')_{(1, 2)}g^1_2G^2)x^1\cd c^*]\\
&&[\smi (f^1\mf{q}^1_1x^3_{(1, 1)}(h_2h')_{(1, 1)}g^1_1G^1)x^2h_1\cd d^*]
\# \smi (\mf{q}^2x^3_2(h_2h')_2g^2)\\
&{{\rm (\ref{tql}, \ref{q1})}\atop =}&
[\smi (\mf{q}^1\mf{Q}^2_1(h_2h')_{(2, 1)}X^2g^1_2G^2)\cd c^*]
[\smi (\mf{Q}^1(h_2h')_1X^1g^1_1G^1)h_1\cd d^*]\\
&&\# \smi (\mf{q}^2\mf{Q}^2_2(h_2h')_{(2, 2)}X^3g^2)\\
\end{eqnarray*}
\begin{eqnarray*}
&{{\rm (\ref{tpql1a}, \ref{g2}, \ref{pf})}\atop =}&
[X^2\smi (\mf{q}^1h_1\mf{Q}^2_1h'_{(2, 1)}g^2_1G^1)\cd c^*]
[X^3\smi (\mf{Q}^1h'_1g^1)\cd d^*]\\
&&\# X^1\smi (\mf{q}^2h_2\mf{Q}^2_2h'_{(2, 2)}g^2_2G^2)\\
&{{\rm (\ref{ca})}\atop =}&
[X^2\smi (\mf{Q}^2h'_2g^2)_2\smi (\mf{q}^1h_1G^1)\cd c^*]
[X^3\smi (\mf{Q}^1h'_1g^1)\cd d^*]\\
&&\# X^1\smi (\mf{Q}^2h'_2g^2)_1\smi (\mf{q}^2h_2G^2)\\
&{{\rm (\ref{stgsm})}\atop =}&\v (c^*\# h)\odot \v (d^*\# h'), 
\end{eqnarray*}
as needed. It is clear that $\v (\une \# 1_H)=\une \# 1_H$, so it
remains to be shown that 
$\v $ is bijective.\\
First we introduce some notation. Let $\mf {A}$ be a right $H$-comodule algebra,
and define the element $\tilde{q}_{\r }\in {\mf A}\ot H$ 
as follows:
\begin{equation}\label{tpqr}
\tilde {q}_{\r }=\tilde {q}^1_{\r }\ot \tilde {q}^2_{\r}
=\tX ^1_{\r }\ot \smi (\a \tX ^3_{\r })\tX ^2_{\r }. 
\end{equation}
In the special situation where
$\mf{A}=H$, the element $\tilde{q}_{\Delta }$ will be denoted by 
$\tilde{q}_{\Delta}=q^1\ot q^2$.\\
We now claim that 
$\varphi ^{-1}:\ ((C^{\rm cop})^*\# \un{H})^{\rm op}\ra C^*\# H$ is given by
the formula
$$
\varphi ^{-1}(c^*\# h)=g^1S(q^2h_2)\cd c^*\# g^2S(q^1h_1), 
$$
for all $c^*\in C^*$ and $h\in H$.\\
From \cite[Lemma 2.6]{bn}, we recall the following formula:
\begin{equation}\label{tfg}
g^2\a \smi (g^1)=\smi (\b ).
\end{equation} 
Let $\tilde{p}_{\Delta }=\mf{p}^1\ot \mf{p}^2$ 
be the element $\tilde{p}_{\l}$ defined in (\ref{tpql}), in the special case where 
$\mf{B}=H$. We then compute, for all $c^*\in C^*$ and $h\in H$:
\begin{eqnarray*}
&&\hspace*{-2.5cm}
\varphi ^{-1}\circ \v (c^*\# h)
=\v ^{-1}(\smi (\mf{q}^1h_1g^1)\cd c^*\# \smi (\mf{q}^2h_2g^2))\\
&{{\rm (\ref{ca})}\atop =}&\mf{q}^2_1h_{(2, 1)}g^2_1G^1S(q^2)
\smi (\mf{q}^1h_1g^1)\cd c^*\# \mf{q}^2_2h_{(2, 2)}g^2_2G^2S(q^1)\\
&{{\rm (\ref{tpqr}, \ref{g2}, \ref{pf})}\atop =}&
\mf{q}^2_1h_{(2, 1)}X^2g^1_2G^2\a \smi (\mf{q}^1h_1X^1g^1_1G^1)\cd c^*
\# \mf{q}^2_2h_{(2, 2)}X^3g^2\\
&{{\rm (\ref{tfg}, \ref{q5}, \ref{tpql})}\atop =}&
\smi (\mf{q}^1h_1S(\mf{q}^2_1h_{(2, 1)}\mf{p}^1))\cd c^*\# 
\mf{q}^2_2h_{(2, 2)}\mf{p}^2\\
&{{\rm (\ref{tpql1}, \ref{tpql2})}\atop =}&c^*\# h.
\end{eqnarray*}
The proof of the fact that 
$\v \circ \v ^{-1}=id_{((C^{\rm cop})^*\# \un{H})^{\rm op}}$ is based 
on similar computations. 
\end{proof}

\begin{remark}\relabel{2.10}
The isomorphism $\v $ in \prref{2.9} can be defined more generally 
for a left $H$-module algebra $A$ instead of a right $H$-module coalgebra 
$C$. Observe that $H$ cannot be replaced by an $H$-bicomodule algebra 
$\mb{A}$, because of the appearance of the antipode $S$ of $H$ on the second 
position of the tensor product. 
\end{remark}

%%%%%%%%%%%%%%%%%%%%%%%%%%%%%%%%%%%%%%
\section{Yetter-Drinfeld modules are Doi-Hopf modules}\selabel{3}
%%%%%%%%%%%%%%%%%%%%%%%%%%%%%%%%%%%%
\setcounter{equation}{0}
In this Section, we will show that 
Yetter-Drinfeld modules are special case of Doi-Hopf modules. 
We will then apply the properties of Doi-Hopf modules 
to Yetter-Drinfeld modules. 

\subsection{Yetter-Drinfeld modules over quasi-bialgebras}\selabel{3.1}
The category of Yetter-Drinfeld modules over a quasi-Hopf algebra $H$ was
 introduced by Majid, as the center of the monoidal category 
${}_H{\cal M}$. His aim was to define the quantum double by an implicit 
Tannaka-Krein reconstruction procedure, see \cite{m1}. From \cite{bpv3},
we recall the following more general definition of Yetter-Drinfeld modules.

The category of $(H,H)$-bimodules, ${}_H{\cal M}_H$, is monoidal. The
associativity constraints ${\bf a'}_{U, V, W}:
(U\ot V)\ot W\ra U\ot (V\ot W)$ are given by 
\begin{equation}\label{bim}
{\bf a'}_{U, V, W}((u\ot v)\ot w)= \Phi \cd (u\ot (v\ot w))\cd\Phi ^{-1}
\end{equation}
for all $U, V, W\in {}_H{\cal M}_H$, $u\in U$, $v\in V$ and $w\in
W$. A coalgebra in the category of
$(H,H)$-bimodules will be called an $H$-bimodule
coalgebra. More precisely, an $H$-bimodule coalgebra $C$ is an
$(H,H)$-bimodule (denote the actions by $h\cd c$ and $c\cd h$)
with a comultiplication $\und :\ C\ra C\ot C$ and
a counit $\une :\ C\ra k$ satisfying the following relations,
for all $c\in C$ and $h\in H$:
\begin{eqnarray}
&&\Phi \cd (\und \ot id_C)(\und (c))\cd \Phi ^{-1}
=(id_C\ot \und )(\und (c)),\label{bmc1}\\
&&\und (h\cd c)=h_1\cd \una \ot h_2\cd \unb ,~~
\und (c\cd h)=\una \cd h_1\ot \unb \cd h_2,\label{bmc2}\\
&&(\une \ot id_C)\circ \und =(id_C\ot \une )\circ \und =id_C,\\
&&\une (h\cd c)=\va (h)\une (c), ~~ 
\une (c\cd h)=\une (c)\va (h),\label{bmc3}
\end{eqnarray}
where we used the same Sweedler-type notation as introduced before.

For further use we note that an $H$-bimodule coalgebra $C$ can be always  
viewed as a left $H^{\rm op}\ot H$-module 
coalgebra via the left $H^{\rm op}\ot H$-action given for all 
$c\in C$ and $h, h'\in H$ by 
\begin{equation}\label{rmcs}
(h\ot h')\cd c=h'\cd c\cd h.
\end{equation}

\begin{definition}(\cite{bpv3})\delabel{3.1}
Let $H$ be a quasi-bialgebra, $C$ an $H$-bimodule coalgebra and 
$\mb {A}$ an $H$-bicomodule algebra. A left-right $(H, \mb{A}, C)$-Yetter-Drinfeld 
module is a $k$-vector space $M$ with the following additional structure:
\begin{itemize}
\item[-] $M$ is a left $\mb {A}$-module; we write $\cdot $ for the 
left $\mb {A}$-action;
\item[-] we have a $k$-linear map $\r _M: M\ra M\ot C$,
$\r _M(m)=m_{(0)}\ot m_{(1)}$, called the right $C$-coaction
on $M$, such that for all $m\in M$, $\une (m_{(1)})m_{(0)}=m$ and
\begin{eqnarray}
&&(\theta ^2\cd m_{(0)})_{(0)}\ot (\theta ^2\cd m_{(0)})_{(1)}\cdot 
\theta ^1\ot \theta ^3\cdot m_{(1)}\nonumber\\
&&\hspace*{1cm}
=\tx ^1_{\r }\cd (\tx ^3_{\l }\cd m)_{(0)}\ot \tx ^2_{\r }\cdot 
(\tx ^3_{\l }\cd m)_{(1)_{\un {1}}}\cdot \tx ^1_{\l }\ot  
\tx ^3_{\r}\cdot (\tx ^3_{\l }\cd m)_{(1)_{\un {2}}}\cdot \tx ^2_{\l },\label{lryd1} 
\end{eqnarray}
\item[-] for all $u\in \mb{A}$ and $m\in M$ 
the following compatibility relation holds:
\begin{eqnarray}
&&u_{<0>}\cd m_{(0)}\ot u_{<1>}\cd m_{(1)}=(u_{[0]}\cd m)_{(0)}
\ot (u_{[0]}\cd m)_{(1)}\cd u_{[-1]}.\label{lryd2}
\end{eqnarray}
\end{itemize}
$_{\mb {A}}{\cal YD}(H)^C$ will be the category of left-right 
$(H, \mb{A}, C)$-Yetter-Drinfeld modules and maps preserving the $\mb {A}$-action 
$C$-coaction.
\end{definition}

We have seen in \seref{1.2} that $H$ is an $H$-bicomodule algebra;
it is also clear that $H$ is an $H$-bimodule coalgebra (take $\ul{\Delta}=\Delta$
and $\ul{\varepsilon}=\varepsilon$). If we take $H=\mb{A}=C$ in
\deref{3.1}, then we recover the category of Yetter-Drinfeld modules
introduced by Majid in \cite{m1}, and studied in \cite{bcp1,bcp2}.
It is remarkable that a quasi-bialgebra $H$ is a coalgebra in the category
of $H$-bimodules, but not in the category of vector spaces, or in the
category of left (or right) $H$-modules.

Let $H$ be a quasi-bialgebra, $\mb{A}$ 
an $H$-bicomodule algebra and $C$ a left $H$-module coalgebra. It is 
straightforward to check that $C$ with the right $H$-module structure 
given by $\va $ is an $H$-bimodule coalgebra. Then 
(\ref{lryd1}) and (\ref{lryd2}) reduce respectively to (\ref{lrdhm1}) and 
(\ref{lrdhm3}), in which only the right $H$-coaction on $\mb{A}$ appears.
So in this 
particular case the category ${}_{\mb{A}}{\cal M}(H)^C$ is just 
${}_{\mb{A}}{\cal YD}(H)^C$.    

In order to show that the Yetter-Drinfeld modules are special case of 
Doi-Hopf modules we need a Doi-Hopf datum. As we have already seen,
an $H$-bimodule coalgebra $C$ can be viewed as a left 
$H^{\rm op}\ot H$-module via the structure defined in (\ref{rmcs}). In order 
to provide a right $H^{\rm op}\ot H$-comodule algebra structure on an $H$-bicomodule 
algebra $\mb{A}$, we need the following result.

\begin{lemma}\lelabel{3.2}
Let $H$ be a quasi-Hopf algebra and $(\mf{B}, \l , \Phi _{\l})$ 
a left $H$-comodule algebra. Then $\mf{B}$ is a right $H^{\rm op}$-comodule algebra 
via the structure
\begin{eqnarray}
&&\r :\ \mf{B}\ra \mf{B}\ot H,~~\r (\mfb )=\mfb _{[0]}\ot \smi (\mfb _{[-1]}),\label{rcastr1}\\
&&\Phi _{\r}=\tilde{x}^3_{\l}\ot \smi (f^2\tilde{x}^2_{\l})\ot \smi (f^1\tilde{x}^1_{\l})
\in \mf{B}\ot H\ot H,\label{rcastr2} 
\end{eqnarray}
where $f^1\ot f^2$ is the Drinfeld twist defined in (\ref{f}). Moreover, if 
$(\mf{B}, \l , \Phi _{\l})$ and $(\mf{B}, \l ', \Phi _{\l '})$ are twist equivalent 
left $H$-comodule algebras then the corresponding right $H^{\rm op}$-comodule algebras 
are also twist equivalent. 
\end{lemma}

\begin{proof}
The relation (\ref{rca1}) follows easily by applying (\ref{lca1}) and (\ref{ca}), and 
the relations (\ref{rca3}, \ref{rca4}) are trivial. We prove now (\ref{rca2}). 
Since $\Phi _{\rm op}=\Phi ^{-1}$ we have 
\begin{eqnarray*}
&&\hspace*{-2cm}
\tilde{X}^1_{\r}\tilde{Y}^1_{\r}\ot x^1\cd _{\rm op}
(\tilde{X}^2_{\r})_1\cd _{\rm op}\tilde{Y}^3_{\r}
\ot x^2\cd _{\rm op}(\tilde{X}^2_{\r})_2
\cd _{\rm op}\tilde{Y}^3_{\r}\ot x^3\cd _{\rm op}\tilde{X}^3_{\r}\\
&=&\tilde{x}^3_{\l}\tilde{y}^3_{\l}\ot 
\smi (F^2\tilde{y}^2_{\l})\smi (f^2\tilde{x}^2_{\l})_1x^1\\
&&\hspace*{1cm}
\ot \smi (F^1\tilde{y}^1_{\l})\smi (f^2\tilde{x}^2_{\l})_2x^2\ot 
\smi (f^1\tilde{x}^1_{\l})x^3\\
&{{\rm (\ref{ca})}\atop =}&
\tilde{x}^3_{\l}\tilde{y}^3_{\l}\ot 
\smi (S(x^1)F^2f^2_2(\tilde{x}^2_{\l})_2\tilde{y}^2_{\l})\\
&&\hspace*{1cm}
\ot \smi (S(x^2)F^1f^2_1(\tilde{x}^2_{\l})_1\tilde{y}^1_{\l})
\ot \smi (S(x^3)f^1\tilde{x}^1_{\l})\\
&{{\rm (\ref{g2}, \ref{pf})}\atop =}&
\tilde{x}^3_{\l}\tilde{y}^3_{\l}\ot 
\smi (F^2x^3(\tilde{x}^2_{\l})_2\tilde{y}^2_{\l})\\
&&\hspace*{1cm}
\ot \smi (f^2F^1_2x^2(\tilde{x}^2_{\l})_1\tilde{y}^1_{\l})
\ot \smi (f^1F^1_1x^1\tilde{x}^1_{\l})\\
&{{\rm (\ref{lca2}, \ref{ca})}\atop =}&
\tilde{y}^3_{\l}(\tilde{x}^3_{\l})_{[0]}\ot 
\smi (F^2\tilde{y}^2_{\l})\cd _{\rm op}
\smi ((\tilde{x}^3_{\l})_{[-1]})\\
&&\hspace*{1cm}
\ot \smi (F^1\tilde{y}^1_{\l})_1\cd _{\rm op}\smi (f^2\tilde{x}^2_{\l})\ot 
\smi (F^1\tilde{y}^1_{\l})_2\cd _{\rm op}\smi (f^1\tilde{x}^1_{\l})\\
&{{\rm (\ref{rcastr1}, \ref{rcastr2})}\atop =}&
\tilde{X}^1_{\r}(\tilde{Y}^1_{\r})_{\le0\ri}\ot \tilde{X}^2_{\r}\cd _{\rm op}
(\tilde{Y}^1_{\r})_{\le1\ri}\ot (\tilde{X}^3_{\r})_1\cd _{\rm op}\tilde{Y}^2_{\r}
\ot (\tilde{X}^3_{\r})_2\cd _{\rm op}\tilde{Y}^3_{\r},
\end{eqnarray*}  
where we denote by $\cd _{\rm op}$ the multiplication of $H^{\rm op}$ and by 
$F^1\ot F^2$ another copy of the Drinfeld twist $f$. 

Finally, it is not hard to see that if 
the invertible element $\mb{U}=\mb{U}^1\ot \mb{U}^2\in H\ot \mf{B}$ 
provides a twist equivalence between the left $H$-comodule 
algebras $(\mf{B}, \l , \Phi _{\l})$ and 
$(\mf{B}, \l ', \Phi _{\l '})$ then the invertible element 
$\mb{V}=\mb{U}^2\ot \smi (\mb{U}^1)\in \mf{B}\ot H$ provides a twist equivalence 
between the associated right $H^{\rm op}$-comodule algebras 
$(\mf{B}, \r , \Phi _{\r})$ and $(\mf{B}, \r ', \Phi _{\r '})$, respectively.   
\end{proof}

\begin{proposition}\prlabel{3.3}
Let $H$ be a quasi-Hopf algebra and $(\mb{A}, \l , \r , \Phi _{\l}, 
\Phi _{\r}, \Phi _{\l , \r})$ an $H$-bicomodule algebra. We define 
two right $H^{\rm op}\ot H$-coactions 
$$
\r _1, \r _2:\ \mb{A}\ra \mb{A}\ot (H^{\rm op}\ot H)
$$ 
on $\mb{A}$, and 
corresponding elements $\Phi _{\r _1}, \Phi _{\r _2}\in 
\mb{A}\ot (H^{\rm op}\ot H)\ot (H^{\rm op}\ot H)$ 
as follows:
\begin{eqnarray}
&&\hspace*{-1cm}
\r _1(u)=u_{\le0\ri _{[0]}}\ot 
\left(\smi (u_{\le0\ri _{[-1]}})\ot u_{\le1\ri}\right),\label{rcaspr1}\\
&&\hspace*{-1cm}
\Phi _{\r _1}=(\tilde{X}^1_{\r})_{[0]}\tilde{x}^3_{\l}\theta ^2_{[0]}\ot 
\left(\smi (f^2(\tilde{X}^1_{\r})_{[-1]_2}\tilde{x}^2_{\l}\theta ^2_{[-1]})\ot 
\tilde{X}^2_{\r}\theta ^3\right) \nonumber\\
&&\hspace*{2.5cm}
\ot \left(\smi (f^1(\tilde{X}^1_{\r})_{[-1]_1}\tilde{x}^1_{\l}\theta ^1)
\ot \tilde{X}^3_{\r}\right),\label{rcaspr2}
\end{eqnarray} 
and 
\begin{eqnarray}
&&\hspace*{-1cm}
\r _2(u)=u_{[0]_{\le0\ri}}\ot \left(\smi (u_{[-1]})
\ot u_{[0]_{\le1\ri}}\right),\label{rcaspr3}\\
&&\hspace*{-1cm}
\Phi _{\r _2}=(\tilde{x}^3_{\l})_{\le0\ri}\tilde{X}^1_{\r}\Theta ^2_{\le0\ri}\ot 
\left(\smi (f^2\tilde{x}^2_{\l}\Theta ^1)\ot 
(\tilde{x}^3_{\l})_{\le1\ri _1}\tilde{X}^2_{\r}
\Theta ^2_{\le1\ri}\right)\nonumber \\
&&\hspace*{2.5cm}
\ot \left(\smi (f^1\tilde{x}^1_{\l})\ot 
(\tilde{x}^3_{\l})_{\le1\ri _2}
\tilde{X}^3_{\r }\Theta ^3\right).\label{rcaspr4}
\end{eqnarray}
Then $(\mb{A}, \r _1, \Phi _{\r _1})$ and 
$(\mb{A}, \r _2, \Phi _{\r _2})$ are twist equivalent 
right $H^{\rm op}\ot H$-comodule algebras. 
\end{proposition}

\begin{proof}
The statement follows from \leref{3.2}. Indeed, we have seen at the end of \seref{1} that 
$\mb{A}$ has two twist equivalent left $H\ot H^{\rm op}$-comodule algebra structures. 
Identifying $(H\ot H^{\rm op})^{\rm op}$ and $H^{\rm op}\ot H$, and computing the induced right 
coactions we obtain the structures defined in (\ref{rcaspr1}-\ref{rcaspr4}). 
We point out that the reassociator, the antipode and the 
Drinfeld twist corresponding to $H\ot H^{\rm op}$ are given by 
\begin{eqnarray*}
&&\Phi _{H\ot H^{\rm op}}=(X^1\ot x^1)\ot (X^2\ot x^2)\ot (X^3\ot x^3),\\  
&&S_{H\ot H^{\rm op}}=S\ot \smi ,~~
f _{H\ot H^{\rm op}}=(f^1\ot \smi (g^2))\ot (f^2\ot \smi (g^1)),
\end{eqnarray*}
where, as usual, $g^1\ot g^2$ is the inverse of $f=f^1\ot f^2$. 
\end{proof}

Let $H$ be a quasi-Hopf algebra, $C$ an $H$-bimodule coalgebra and $\mb{A}$ an 
$H$-bicomodule algebra. In the sequel, $\mb{A}^1$ and 
$\mb{A}^2$ will be our notation for the right $H^{\rm op}\ot H$-comodule algebras 
$(\mb{A}, \r _1, \Phi _{\r _1})$ 
and $(\mb{A}, \r _2, \Phi _{\r _2})$. By the above arguments, it make sense 
to consider the left-right Doi-Hopf module categories 
${}_{\mb{A}^1}{\cal M}(H^{\rm op}\ot H)^C$ and ${}_{\mb{A}^2}{\cal M}(H^{\rm op}\ot H)^C$.
It will follow from \prref{3.4} that these two categories are isomorphic. 

\begin{proposition}\prlabel{3.4}
Let $H$ be a quasi-bialgebra, $C$ a left $H$-module coalgebra 
and $\mfA ^1=(\mfA , \r , \Phi _{\r})$ and 
$\mfA ^2=(\mfA , \r ', \Phi _{\r '})$ two twist equivalent right $H$-comodule algebras.
Then the categories ${}_{\mfA ^1}{\cal M}(H)^C$ and 
${}_{\mfA ^2}{\cal M}(H)^C$ are isomorphic.
\end{proposition} 

\begin{proof}
If $\mfA ^1$ and $\mfA ^2$ are twist equivalent, then there exists
$\mb{V}\in \mfA \ot H$ satisfying (\ref{eq:comtwist0}-\ref{eq:comtwist2}).
Take $M\in {}_{\mfA ^1}{\cal M}(H)^C$; $M$ becomes 
an object in  ${}_{\mfA ^2}{\cal M}(H)^C$ by keeping the same left $\mfA$-module 
structure and defining
$$
\r '_M:\ M\ra M\ot C,~~\r '(m)=\mb{V}\cd \r _M(m).$$
Conversely, take $M\in {}_{\mfA ^2}{\cal M}(H)^C$ via the structures 
$\cd $ and $\r '_M$. Then $M$ can be viewed as a left-right $(H, \mfA ^1, C)$-Hopf 
module
via the same left $\mfA$-action $\cd $ and the right $C$-coaction $\r _M$ defined by   
$$
\r _M:\ M\ra M\ot C,~~\r _M(m)=\mb{V}^{-1}\cd \r '_M(m).
$$
These correspondences define two functors which act as the identity on morphisms and 
produce inverse isomorphisms.
\end{proof}

\begin{remark}\relabel{3.5}
Let $H$ be a quasi-bialgebra, and $F\in  H\ot H$ a gauge transformation.
We can consider the twisted quasi-bialgebra $H_F$ (see (\ref{g1}-\ref{g2})), 
and it is well-known that the categories of left $H$-modules and left $H_F$-modules 
are isomorphic. We have a similar property for Doi-Hopf modules.\\
Let $(\mfA, \r , \Phi _{\r})$ be a right $H$-comodule algebra, and 
let $\Phi _{\r _F}=(1_{\mfA}\ot F)\Phi _{\r}$. Then 
$\mfA _F:=(\mfA , \r , \Phi _{\r _F})$ is a right 
$H_F$-comodule algebra, see \cite{hn1}.\\
Now let $C$ be a left $H$-module coalgebra, and define a new
comultiplication $\und _{F}$ as follows: $\und _{F}(c)=F\und (c)$,
for all $c\in C$. Then straightforward computations show that
$(C, \und _{F}, \une )$ is a left $H_F$-module coalgebra and
that the categories 
${}_{\mfA}{\cal M}(H)^C$ and ${}_{\mfA _F}{\cal M}(H_F)^{C_F}$ are isomorphic. Of 
course, a similar result holds for left comodule algebras. 
\end{remark}

\subsection{Yetter-Drinfeld modules and Doi-Hopf modules}\selabel{3.2}
Our next aim is to show that the category of left-right Yetter-Drinfeld 
modules ${}_{\mb{A}}{\cal YD}(H)^C$ 
is isomorphic to the category of Doi-Hopf modules 
${}_{\mb{A}^2}{\cal M}(H^{\rm op}\ot H)^C$ and, a fortiori,  
to ${}_{\mb{A}^1}{\cal M}(H^{\rm op}\ot H)^C$, by \prref{3.4}. 
We have divided the proof over a few lemmas.

\begin{lemma}\lelabel{3.6}
Let $H$ be a quasi-Hopf algebra, $\mb{A}$ an $H$-bicomodule algebra and 
$C$ an $H$-bimodule coalgebra. We have a functor
$$
F:\ {}_{\mb{A}}{\cal YD}(H)^C\ra {}_{\mb{A}^2}{\cal M}(H^{\rm op}\ot H)^C
$$
which acts as the identity on objects an morphisms. If $M$ is a left-right 
$(H, \mb{A}, C)$-Yetter-Drinfeld module then $F(M)=M$ as a left 
$\mb{A}$-module, and with the newly defined right $C$-coaction 
\begin{equation}\label{ydstru1}
\r '_M(m)=m_{(0)'}\ot m_{(1)'}=(\tplb \cd m)_{(0)}\ot (\tplb \cd m)_{(1)}\cd \tpla,
\end{equation} 
for all $m\in M$. Here 
$\tilde{p}_{\l}=\tpla \ot \tplb$ is the element defined by (\ref{tpql}). 
\end{lemma}

\begin{proof}
It is not hard to see that (\ref{bca2}) and (\ref{q5}) imply
\begin{equation}\label{fs1}
\Theta ^2_{[0]}\tplb \ot \Theta ^2_{[-1]}\tpla \smi (\Theta ^1)\ot \Theta ^3
=\theta ^2(\tplb )_{\le0\ri}\ot \theta ^1\tpla \ot \theta ^3(\tplb )_{\le1\ri}.
\end{equation}
We now show that $F(M)$ satisfies the relations (\ref{lrdhm1}-\ref{lrdhm3}). Let 
$\tPla \ot \tPlb$ be another copy of $\tilde{p}_{\l}$, and compute
\begin{eqnarray*}
&&\hspace*{-2.5cm}
\Phi _{\r _2}\cd (\r '_M\ot id_C)(\r _M(m))\\
&{{\rm (\ref{rcaspr4}, \ref{ydstru1}, \ref{rmcs})}\atop =}&
(\tx ^3_{\l})_{\le0\ri}\tX ^1_{\r}\Theta ^2_{\le0\ri}\cd 
(\tPlb \cd (\tplb \cd m)_{(0)})_{(0)}\\
&&\ot (\tx ^3_{\l})_{\le1\ri _1}\tX ^2_{\r}\Theta ^2_{\le1\ri}\cd 
(\tPlb \cd (\tplb \cd m)_{(0)})_{(1)}\cd 
\tPla \smi (f^2\tx ^2_{\l}\Theta ^1)\\
&&\ot (\tx ^3_{\l})_{\le1\ri _2}\tX ^3_{\r}\Theta ^3\cd (\tplb \cd m)_{(1)}\cd 
\tpla \smi (f^1\tx ^1_{\l})\\
&{{\rm (\ref{lryd2}, \ref{fs1})}\atop =}&
(\tx ^3_{\l})_{\le0\ri}\tX ^1_{\r}\cd (\theta ^2(\tPlb )_{\le0\ri}\cd 
(\tplb \cd m)_{(0)})_{(0)}\\
&&\ot (\tx ^3_{\l})_{\le1\ri _1}\tX ^2_{\r}\cd (\theta ^2(\tPlb )_{\le0\ri}\cd 
(\tplb \cd m)_{(0)})_{(1)}\cd \theta ^1\tPla \smi (f^2\tx ^2_{\l})\\
&&\ot (\tx ^3_{\l})_{\le1\ri _2}\tX ^3_{\r}\theta ^3(\tPlb )_{\le1\ri}\cd 
(\tplb \cd m)_{(1)}\cd \tpla \smi (f^1\tx ^1_{\l})\end{eqnarray*}
\begin{eqnarray*}
&{{\rm (\ref{lryd2}, \ref{lryd1})}\atop =}&
(\tx ^3_{\l})_{\le0\ri}\cd (\ty ^3_{\l}(\tPlb )_{[0]}\tplb \cd m)_{(0)}\\
&&\ot (\tx ^3_{\l})_{\le1\ri _1}\cd (\ty ^3_{\l}(\tPlb )_{[0]}\tplb \cd m)_{(1)_{\un{1}}}
\cd \ty ^1_{\l}\tPla\smi (f^2\tx ^2_{\l})\\
&&\ot (\tx ^3_{\l})_{\le1\ri _2}\cd (\ty ^3_{\l}(\tPlb )_{[0]}\tplb \cd m)_{(1)_{\un{2}}}\cd 
\ty ^2_{\l}(\tPlb )_{[-1]}\tpla \smi (f^1\tx ^1_{\l})\\
&{{\rm (\ref{bmc2}, \ref{lryd2}, \ref{tpl})}\atop =}&
(\tplb \cd m)_{(0)}\ot (\tplb \cd m)_{(1)_{\un{1}}}\cd (\tpla )_1\ot 
(\tplb \cd m)_{(1)_{\un{2}}}\cd (\tpla )_2\\
&{{\rm (\ref{bmc2}, \ref{ydstru1})}\atop =}&
m_{(0)'}\ot m_{(1)'_{\un{1}}}\ot m_{(1)'_{\un{2}}}=(id_M\ot \und )(\r '_M(m)),
\end{eqnarray*}   
for all $m\in M$, as needed. The relation (\ref{lrdhm2}) is trivial and 
(\ref{lrdhm3}) follows from (\ref{tpql1}), (\ref{rcaspr3}) and (\ref{rmcs}).  
\end{proof}

\begin{lemma}\lelabel{3.7}
Let $H$ be a quasi-Hopf algebra, $C$ an $H$-bimodule coalgebra and 
$\mb{A}$ an $H$-bicomodule algebra. Then we have a functor 
$$
G:\ {}_{\mb{A}^2}{\cal M}(H^{\rm op}\ot H)^C\ra {}_{\mb{A}}{\cal YD}(H)^C
$$
which acts as the identity on objects and morphisms. Let $M$ be a left-right 
$(H^{\rm op}\ot H, \mb{A}^2, C)$-Hopf module, with  left 
$\mb{A}$-action $\cd$ and right $C$-coaction $\rho'_M$,
$\r '_M(m)=m_{(0)'}\ot m_{(1)'}\in M\ot C$. Then $G(M)=M$ 
as a left $\mb{A}$-module, with new right $C$-coaction 
$\ov{\r}_M:\ M\ra M\ot C$, given by the formula
\begin{equation}\label{ydstrsp2}
\ov{\r}_M(m)=m_{\ov{(0)}}\ot m_{\ov{(1)}}=
(\tqlb )_{\le0\ri}\cd m_{(0)'}\ot (\tqlb )_{\le1\ri}\cd m_{(1)'}\cd \smi (\tqla ),
\end{equation} 
for all $m\in M$. Here
$\tilde{q}_{\l}=\tqla \ot \tqlb$ is the element defined in (\ref{tpql}).
\end{lemma}

\begin{proof}
The most difficult part is tho show that $G(M)$ satisfies the relations 
(\ref{lryd1}) and (\ref{lryd2}). $M$ is a left-right $(H, \mb{A}^2, C)$-Hopf 
module, so we have
 by (\ref{lrdhm1}), (\ref{rcaspr3}, \ref{rcaspr4}) and (\ref{rca1}):
\begin{eqnarray}
&&\tX ^1_{\r}((\tx ^3_{\l})_{\le0\ri}\Theta ^2)_{\le0\ri}\cd m_{(0, 0)'}\ot 
\tX ^2_{\r}((\tx ^3_{\l})_{\le0\ri}\Theta ^2)_{\le1\ri}\cd m_{(0, 1)'}
\cd \smi (f^2\tx ^2_{\l}\Theta ^1)\nonumber\\
&&\hspace*{2cm}
\ot \tX ^3_{\r}(\tx ^3_{\l})_{\le1\ri}\Theta ^3\cd m_{(1)'}\cd 
\smi (f^1\tx ^1_{\l})=m_{(0)'}\ot m_{(1)'_{\un{1}}}\ot m_{(1)'_{\un{2}}},\label{fs2}\\
&&(u\cd m)_{(0)'}\ot (u\cd m)_{(1)'}=u_{[0]_{\le0\ri}}\cd m_{(0)'}\ot 
u_{[0]_{\le1\ri}}\cd m_{(1)'}\cd \smi (u_{[-1]}),\label{fs3}
\end{eqnarray}
for all $m\in M$ and $u\in \mb{A}$. Also, (\ref{bca2}) and (\ref{q5}) imply that
\begin{equation}\label{fs4}
S(\theta ^1)\tqla \theta ^2_{[-1]}\ot \tqlb \theta ^2_{[0]}\ot \theta ^3=
\tqla \Theta ^1\ot (\tqlb )_{\le0\ri}\Theta ^2\ot (\tqlb )_{\le1\ri}\Theta ^3.
\end{equation}

Let $\tQla \ot \tQlb$ be another copy of 
$\tilde{q}_{\l}$; for all $m\in M$, we compute that
\begin{eqnarray*}
&&\hspace*{-2.5cm}
(\theta ^2\cd m_{\ov{(0)}})_{\ov{(0)}}\ot (\theta ^2\cd m_{\ov{(0)}})_{\ov{(1)}}
\ot \theta ^3\cd m_{\ov{(1)}}\\
&{{\rm (\ref{ydstrsp2}, \ref{fs3})}\atop =}&
(\tQlb \theta ^2_{[0]}(\tqlb )_{\le0\ri _{[0]}})_{\le0\ri}\cd m_{(0, 0)'}
\ot (\tQlb \theta ^2_{[0]}(\tqlb )_{\le0\ri _{[0]}})_{\le1\ri}\cd m_{(0, 1)'}\\
&&\cd \smi (\tQla \theta ^2_{[-1]}(\tqlb )_{\le0\ri _{[-1]}})\theta ^1\ot 
\theta ^3(\tqlb )_{\le1\ri}\cd m_{(1)'}\cd \smi (\tqla )\\
&{{\rm (\ref{fs4}, \ref{bca1})}\atop =}&
((\tQlb (\tqlb )_{[0]})_{\le0\ri}\Theta ^2)_{\le0\ri}\cd m_{(0, 0)'}\ot 
((\tQlb (\tqlb )_{[0]})_{\le0\ri}\Theta ^2)_{\le1\ri}\cd m_{(0, 1)'}\\
&&\cd \smi (\tQla (\tqlb )_{[-1]}\Theta ^1)\ot 
(\tQlb (\tqlb )_{[0]})_{\le1\ri}\Theta ^3\cd m_{(1)'}\cd \smi (\tqla )
\end{eqnarray*}
\begin{eqnarray*}
&{{\rm (\ref{tql}, \ref{fs2}, \ref{ca})}\atop =}& 
(\tqlb (\tx ^3_{\l})_{[0]})_{\le0, 0\ri}\tx ^1_{\r}\cd m_{(0)'}\\
&&\ot (\tqlb (\tx ^3_{\l})_{[0]})_{\le0, 1\ri}\tx ^2_{\r}\cd m_{(1)'_{\un{1}}}
\cd \smi (\tqla (\tx ^3_{\l})_{[-1]})_1\tx ^1_{\l}\\
&&\ot (\tqlb (\tx ^3_{\l})_{[0]})_{\le1\ri}\tx ^3_{\r }\cd m_{(1)'_{\un{2}}}\cd 
\smi (\tqla (\tx ^3_{\l})_{[-1]})_2\tx ^2_{\l}\\
&{{\rm (\ref{rca1}, \ref{bca2}, \ref{fs3})}\atop =}&
\tx ^1_{\r}\cd ((\tqlb )_{\le0\ri}\cd (\tx ^3_{\l}\cd m)_{(0)'})\\
&&\ot \tx ^2_{\r}\cd \left((\tqlb )_{\le1\ri}\cd 
(\tx ^3_{\l}\cd m)_{(1)'}\cd \smi (\tqla )\right)_{\un{1}}\cd \tx ^1_{\l}\\
&&\ot \tx ^3_{\r}\cd \left((\tqlb )_{\le1\ri}\cd 
(\tx ^3_{\l}\cd m)_{(1)'}\cd \smi (\tqla )\right)_{\un{2}}\cd \tx ^2_{\l}\\
&{{\rm (\ref{ydstrsp2})}\atop =}&\tx ^1_{\r}\cd (\tx ^3_{\l}\cd m)_{\ov{(0)}}\ot 
\tx ^2_{\r}\cd (\tx ^3_{\l}\cd m)_{\ov{(1)}_{\un{1}}}\cd \tx ^1_{\l}\ot 
\tx ^3_{\r}\cd (\tx ^3_{\l}\cd m)_{\ov{(1)}_{\un{2}}}\cd \tx ^2_{\l}, 
\end{eqnarray*}
as needed. (\ref{lryd2}) also holds since 
\begin{eqnarray*}
&&\hspace*{-2.5cm}
u_{\le0\ri}\cd m_{\ov{(0)}}\ot u_{\le1\ri}\cd m_{\ov{(1)}}\\
&{{\rm (\ref{ydstrsp2}, \ref{tpql1a})}\atop =}&
(\tqlb u_{[0, 0]})_{\le0\ri}\cd m_{(0)'}\ot (\tqlb u_{[0, 0]})_{\le1\ri}
\cd m_{(1)'}\cd \smi (\tqla u_{[0, -1]})u_{[-1]}\\
&{{\rm (\ref{fs3})}\atop =}&
(\tqlb )_{\le0\ri}\cd (u_{[0]}\cd m)_{(0)'}\ot (\tqlb )_{\le1\ri}\cd 
(u_{[0]}\cd m)_{(1)'}\cd \smi (\tqla )u_{[-1]}\\
&{{\rm (\ref{ydstrsp2})}\atop =}&(u_{[0]}\cd m)_{\ov{(0)}}\ot 
(u_{[0]}\cd m)_{\ov{(1)}}\cd u_{[-1]},
\end{eqnarray*}
for all $u\in \mb{A}$ and $m\in M$. The remaining details are left to the reader. 
\end{proof}

\thref{3.8} is the main result of this Section, and
generalizes \cite[Theorem 2.3]{cmz} to the quasi-Hopf algebra setting.

\begin{theorem}\thlabel{3.8}
Let $H$ be a quasi-Hopf algebra, $\mb{A}$ an $H$-bicomodule algebra 
and $C$ an $H$-bimodule coalgebra. Then the categories 
${}_{\mb{A}}{\cal YD}(H)^C$ and ${}_{\mb{A}^2}{\cal M}(H^{\rm op}\ot H)^C$ 
are isomorphic. In particular ${}_{\mb{A}}{\cal YD}(H)^C$ is a Grothendieck 
category, and therefore it has enough injective objects. 
\end{theorem}

\begin{proof}
We show that the functors $F$ and $G$ from Lemmas \ref{le:3.6} and \ref{le:3.7} 
are inverses. Since $F$ and $G$ act as the identity functor at the level 
of $\mb{A}$-modules we have only to show that $F$ and $G$ are inverses at 
the level of $C$-coactions. 

Let $M\in {}_{\mb{A}}{\cal YD}(H)^C$ and $\r _M(m)=m_{(0)}\ot m_{(1)}$ 
its right $C$-coaction. We denote by $\ov{\r}_M(m)=m_{\ov{(0)}}\ot m_{\ov{(1)}}$ 
the right $C$-coaction of $G(F(M))$ obtained using first \leref{3.6} 
and then \leref{3.7}. For all $m\in M$ we then have
\begin{eqnarray*}
\ov{\r}_M(m)
&{{\rm (\ref{ydstrsp2})}\atop =}&
(\tqlb )_{\le0\ri}\cd m_{(0)'}\ot (\tqlb )_{\le1\ri}\cd m_{(1)'}\cd \smi (\tqla )\\
&{{\rm (\ref{ydstru1})}\atop =}&
(\tqlb )_{\le0\ri}\cd (\tplb \cd m)_{(0)}\ot 
(\tqlb )_{\le1\ri}\cd (\tplb \cd m)_{(1)}\cd \tpla \smi (\tqla )\\
&{{\rm (\ref{lryd2})}\atop =}&
((\tqlb )_{[0]}\tplb \cd m)_{(0)}\ot ((\tqlb )_{[0]}\tplb \cd m)_{(1)}\cd 
(\tqlb )_{[-1]}\tpla \smi (\tqla )\\
&{{\rm (\ref{tpql2})}\atop =}&m_{(0)}\ot m_{(1)}=\r _M(m).
\end{eqnarray*}

Conversely, let $M\in {}_{\mb{A}^2}{\cal M}(H^{\rm op}\ot H)^C$ and denote by 
$\r '_M(m)=m_{(0)'}\ot m_{(1)'}$ its right $C$-coaction. If 
$\r _M(m)=m_{(0)}\ot m_{(1)}$ is the right $C$-coaction on $F(G(M))$ obtained 
using first \leref{3.7} and then \leref{3.6} we have 
\begin{eqnarray*}
\r _M(m)&{{\rm (\ref{ydstru1})}\atop =}&(\tplb \cd m)_{\ov{(0)}}\ot 
(\tplb \cd m)_{\ov{(1)}}\cd \tpla \\
&{{\rm (\ref{ydstrsp2})}\atop =}&(\tqlb )_{\le0\ri}\cd (\tplb \cd m)_{(0)'}\ot 
(\tqlb )_{\le1\ri}\cd (\tplb \cd m)_{(1)'}\cd \smi (\tqla )\tpla \\
&{{\rm (\ref{fs3})}\atop =}&(\tqlb (\tplb )_{[0]})_{\le0\ri}\cd m_{(0)'}\ot 
(\tqlb (\tplb )_{[0]})_{\le1\ri}\cd m_{(1)'}\cd \smi (\tqla (\tplb )_{[-1]})\tpla \\
&{{\rm (\ref{tpql2a})}\atop =}&m_{(0)'}\ot m_{(1)'}=\r '_M(m),
\end{eqnarray*}
for all $m\in M$, so the proof is finished. 
\end{proof}

Let $H$ be a quasi-bialgebra, $\mfA$ a right $H$-comodule algebra and $C$ a left 
$H$-module coalgebra. Identifying the category of left-right $(H, \mfA , C)$-Hopf modules 
to the category of right-left $(H^{\rm op,cop}, \un{\mfA}^{\rm op}, C^{\rm cop})$-Hopf 
modules using the construction preceding \prref{2.2}, we obtain the functor,
after permuting the tensor factors: 
$${\cal F}'=\bullet \ot C:\ {}_{\mfA}{\cal M}\ra {}_{\mfA}{\cal M}(H)^C.$$
 If $M$ is a left 
$\mfA$-module then ${\cal F}'(M)=M\ot C$ with structure  maps
\begin{eqnarray*}
&&\mfa \cd (m\ot c)=\mfa _{\le0\ri}\cd m\ot \mfa _{\le1\ri}\cd c,\\
&&\r _M(m\ot c)=\tx ^1_{\r}\cd m\ot \tx ^2_{\r}\cd \una \ot \tx ^3_{\r}\cd \unb,  
\end{eqnarray*}
for all $\mfa \in \mfA$, $m\in M$ and $c\in C$. For a morphism $\upsilon $
in  ${}_{\mfA}{\cal M}$,
we have that ${\cal F}'(\upsilon )=\upsilon \ot id_C$. In particular, we obtain that 
$\mfA \ot C$ is a left-right $(H, \mfA , C)$-Hopf module.  

Moreover, ${\cal F}'$ is a right adjoint of the forgetful functor 
${\cal U}^C:\ {}_{\mfA}{\cal M}(H)^C\ra {}_{\mfA}{\cal M}$, and it is a left adjoint 
of the functor 
$\Hom_{\mfA}^C(\mfA \ot C, \bullet ):\ {}_{\mfA}{\cal M}(H)^C\ra {}_{\mfA}{\cal M}$ 
defined as follows. For $M\in {}_{\mfA}{\cal M}(H)^C$, 
$\Hom_{\mfA}^C(\mfA \ot C, \bullet )(M)=\Hom_{\mfA}^C(\mfA \ot C, M)$, the set of 
morphisms in ${}_{\mfA}{\cal M}(H)^C$ between $\mfA \ot C$ and $M$, viewed as a 
left $\mfA$-module via 
$$
(\mfa \cd \eta )(\mfa '\ot c)=\eta (\mfa '\mfa \ot c),$$
for all
$\eta \in \Hom_{\mfA}^C(\mfA \ot C, M)$, $\mfa , \mfa '\in \mfA$ and $c\in C$.
The functor $\Hom_{\mfA}^C(\mfA \ot C, \bullet )$ sends a morphism $\kappa$ 
from ${}_{\mfA}{\cal M}(H)^C$ to the morphism 
$\vartheta \mapsto \kappa \circ \vartheta$. 

\begin{corollary}\colabel{3.9}
Let $H$ be a quasi-Hopf algebra, $\mb{A}$ an $H$-bicomodule algebra and $C$ an 
$H$-bimodule coalgebra. We have a functor 
$\mf{F}=\bullet \ot C:\ {}_{\mb{A}}{\cal M}\ra {}_{\mb{A}}{\cal YD}(H)^C$. 
The structure maps on $\mf{F}(M)=M\ot C$ 
are the following:
\begin{eqnarray*}
&&u\cd (m\ot c)=u_{[0]_{\le0\ri}}\cd m\ot u_{[0]_{\le1\ri}}\cd c\cd \smi (u_{[-1]}),\\
&&\r _{M\ot C}(m\ot c)=\theta ^2_{\le0\ri}\tx ^1_{\r}((\tqlb )_{[0]}
\tX ^3_{\l})_{\le0\ri}\cd m
\ot \theta ^2_{\le1\ri}\tx ^2_{\r}((\tqlb )_{[0]}\tX ^3_{\l})_{\le1\ri _1}\cd \una \\
&&\hspace*{1.3cm}
\cd \smi (\theta ^1(\tqlb )_{[-1]}\tX ^2_{\l}g^2)
\ot \theta ^3\tx ^3_{\r}((\tqlb )_{[0]}\tX ^3_{\l})_{\le1\ri _2}\cd \unb \cd 
\smi (\tqla \tX ^1_{\l}g^1),
\end{eqnarray*} 
for all $u\in \mb{A}$, $m\in M$ and $c\in C$. In particular, $\mb{A}\ot C$ 
is a left-right $(H, \mb{A}, C)$-Yetter-Drinfeld module. 
Moreover, the following assertions hold:
\begin{itemize}
\item[i)] $\mf{F}$ is right adjoint to the forgetful functor 
$\mf{U}^C:\ {}_{\mb{A}}{\cal YD}(H)^C\ra {}_{\mb{A}}{\cal M}$.
\item[ii)] $\mf{F}$ is left adjoint to the functor 
$\Hom_{\mb{A}}^C(\mb{A}\ot C, \bullet ):\ {}_{\mb{A}}{\cal YD}(H)^C\ra {}_{\mb{A}}{\cal M}$. 
If $M\in {}_{\mb{A}}{\cal YD}(H)^C$ then 
$\Hom_{\mb{A}}^C(\mb{A}\ot C, \bullet )(M)=\Hom_{\mb{A}}^C(\mb{A}\ot C, M)$, 
the set of Yetter-Drinfeld morphisms from $\mb{A}\ot C$ to $M$, viewed as a 
right $\mb{A}$-module via the action 
\begin{eqnarray*}
&&\hspace*{5mm}
(u \cd \eta )(u'\ot c)=\eta (u'u\ot c),
\end{eqnarray*}
for all
$\eta \in \Hom_{\mb{A}}^C(\mb{A} \ot C, M)$, $u, u'\in \mb{A}$ and $c\in C$.
\end{itemize}
\end{corollary}

\begin{proof}
The functor $\mf{F}$ is well defined because it is the composition 
\begin{diagram}
\mf{F}: {}_{\mb{A}}{\cal M}={}_{\mb{A}^2}{\cal M}&\pile{\rTo^{{\cal F}'}}
&{}_{\mb{A}^2}{\cal M}(H^{\rm op}\ot H)^C&\pile{\rTo^{G}}&{}_{\mb{A}}{\cal YD}(H)^C,
\end{diagram} 
where ${\cal F}'$ is the functor described above but now for the context 
given by (\ref{rcaspr3}-\ref{rcaspr4}), and $G$ is the functor from \leref{3.7}. 
All the others details are left to the reader. 
\end{proof}

Let $H$ be a quasi-bialgebra, $\mfA$ a right $H$-comodule algebra 
and $C$ a left $H$-module coalgebra. Since the category 
${}_{\mfA}{\cal M}(H)^C$ can be identified to a category of right-left 
Doi-Hopf modules it follows that ${}_{\mfA}{\cal M}(H)^C$ is isomorphic to 
a category of comodules over the coring ${\cal C}'=C\ot \mfA$,
with $\mfA$-bimodule structure given by 
$$
\mfa \cd (c\ot \mfa ')\cd \mfa {''}=\mfa _{\le1\ri}\cd c\ot 
\mfa _{\le0\ri}\mfa '\mfa {''}
$$
and comultiplication and counit given by 
$$
\Delta _{{\cal C}'}(c\ot \mfa )=(\tx ^3_{\r}\cd \unb \ot 1_{\mfA})\ot _{\mfA}
(\tx ^2_{\r}\cd \una \ot \tx ^1_{\r}\mfa ),~~
\va _{{\cal C}'}(c\ot \mfa )=\une (c)\mfa 
$$
for all $\mfa , \mfa ', \mfa {''}\in \mfA$ and $c\in C$. Using arguments 
similar to the ones given 
in \cite[Theorem 5.4]{bc3}, one can easily 
check that ${}_{\mfA}{\cal M}(H)^C$ is isomorphic to 
${}^{{\cal C}'}{\cal M}$.

\begin{corollary}\colabel{3.10}
Let $H$ be a quasi-Hopf algebra, $\mb{A}$ an $H$-bicomodule 
algebra and $C$ an $H$-bimodule coalgebra. Then the category of 
left-right Yetter-Drinfeld modules ${}_{\mb{A}}{\cal YD}(H)^C$ 
is isomorphic to the category of left comodules over the 
coring $\mf{C}'=C\ot \mb{A}$, with the following structure maps. 
The $\mb{A}$-bimodule 
is given by 
$$
u\cd (c\ot u')\cd u{''}=u_{[0]_{\le1\ri}}\cd c\cd \smi (u_{[-1]})\ot 
u_{[0]_{\le0\ri}}u'u'',
$$
and the comultiplication and counit are defined by the formulas
\begin{eqnarray*}
&&\hspace*{-1cm}
\Delta _{\mf{C}'}(c\ot u)=\left(\theta ^3\tx ^3_{\r}(\tX ^3_{\l})_{\le1\ri _2}\cd \unb 
\cd \smi (\tX ^1_{\l}g^1)\ot 1_{\mb{A}}\right)\\
&&\ot _{\mb{A}}
\left(\theta ^2_{\le1\ri}\tx ^2_{\r}(\tX ^3_{\l})_{\le1\ri _1}\cd \una \cd 
\smi (\theta ^1\tX ^2_{\l}g^2)
\ot \theta ^2_{\le0\ri}\tx ^1_{\r}(\tX ^3_{\l})_{\le0\ri}u\right),\\
&&\hspace*{-1cm}
\va _{\mf{C}'}(c\ot u)=\une (c)u,
\end{eqnarray*} 
for all $u, u', u''\in \mb{A}$ and $c\in C$.
\end{corollary} 

\begin{proof}
This is a direct consequence of the above comments and \thref{3.8}. 
\end{proof}

\subsection{The category of Yetter-Drinfeld modules as a module category}\selabel{3.3}
Our next aim is to describe the category of Yetter-Drinfeld modules
as a category of modules. We will need the (right) generalized diagonal 
crossed product construction, as introduced in \cite{hn1,bpv3}.

Let $H$ be a quasi-bialgebra and $\mb{A}$ an $H$-bicomodule algebra.
In the sequel, ${\cal A}$ will be an $H$-bimodule algebra. 
This means that ${\cal A}$ is an 
$H$-bimodule which has a multiplication and a usual unit $1_{\cal A}$ 
such that for all $\varphi , \psi , \xi \in {\cal A}$ and 
$h\in H$ the following relations hold:
\begin{eqnarray*}
&&\hspace{-2cm}
(\varphi \psi )\xi =(X^1\cd \varphi \cd x^1)[(X^2\cd \psi \cd x^2)
(X^3\cd \xi \cd x^3)],\label{bma1}\\
&&\hspace{-2cm}
h\cd (\varphi \psi)=(h_1\cd \varphi)(h_2\cd \psi ), ~~
(\varphi \psi )\cd h=(\varphi \cd h_1)(\psi \cd h_2),\label{bma2}\\
&&\hspace{-2cm}
h\cd 1_{\cal A}=\va (h)1_{\cal A},~~
1_{\cal A}\cd h=\va (h)1_{\cal A}.\label{bma3}
\end{eqnarray*}
If $H$ is a quasi-bialgebra, then $H^*$,
the linear dual of $H$, is an $H$-bimodule via the 
$H$-actions 
$$\le h\rh \v , h{'}\ri=\v (h{'}h), ~~
\le \v \lh h, h{'}\ri=\v (hh{'}),$$
for all $\varphi \in H^*$, $h, h{'}\in H$. The convolution 
$\le \v \psi , h\ri=\v (h_1)\psi (h_2)$, $\v ,
\psi \in H^*$, $h\in H$, is a multiplication on $H^*$ which is not
associative in general, but with this multiplication $H^*$ becomes
an $H$-bimodule algebra.

Let $H$ be a quasi-bialgebra, ${\cal A}$ an $H$-bimodule algebra and 
$(\mb{A}, \l , \r , \Phi _{\l}, \Phi _{\r}, \Phi _{\l , \r})$ 
an $H$-bicomodule algebra. 
In the sequel $(\d , \Psi )$ will be the pair 
\begin{eqnarray*}
&&\d _l=(\l \ot id_H)\circ \r ,\label{dl}\\
&&\Psi _:=(id_H\ot \l \ot id_H^{\ot 2})\left((\Phi _{\l ,\r}\ot 1_H)(\l \ot id_H^{\ot 2})
(\Phi _{\r}^{-1})\right)[\Phi _{\l}\ot 1_H^{\ot 2}],\label{pl}
\end{eqnarray*} 
or  
\begin{eqnarray*}
&&\d _r=(id_H\ot \r )\circ \l ,\label{dr}\\
&&\Psi _{r}=(id_H^{\ot 2}\ot \r \ot id_H)\left((1_H\ot \Phi _{\l ,\r}^{-1})
(id_H^{\ot 2}\ot \r )(\Phi _{\l})\right)[1_H^{\ot 2}\ot \Phi _{\r}^{-1}].
\end{eqnarray*}
$\Omega _{L_{l/r}}, \Omega _{R_{l/r}}\in H^{\ot 2}\ot \mb{A}\ot H^{\ot 2}$ 
are defined by the following formulas  
\begin{eqnarray*}
&&\Omega _{L_{l/r}}=(id_H^{\ot 2}\ot id_{\mb{A}}\ot \smi \ot \smi )(\Psi ^{-1}_{l/r})\cd 
[1_H^{\ot 2}\ot 1_{\mb{A}}\ot \smi (f^1)\ot \smi (f^2)],\label{ol}\\ 
&&\Omega _{R_{l/r}}=[\smi (g^1)\ot \smi (g^2)\ot 1_{\mb{A}}\ot 1_H^{\ot 2}]\cd 
(\smi \ot \smi \ot id_{\mb{A}}\ot id_H^{\ot 2})(\Psi _{l/r}).\label{or}
\end{eqnarray*}
Here $f=f^1\ot f^2$ is the Drinfeld twist and $f^{-1}=g^1\ot g^2$ 
is its inverse. We will use the notation
$$\Omega _{L_{l/r}}=\O ^1_{L_{l/r}}\ot \cdots \ot \O ^5_{L_{l/r}}.$$
Let ${\cal A}\bowtie _{\d _l}{\mb A}$ and ${\cal A}\btrl _{\d _r}{\mb A}$
be the vector space ${\cal A}\ot \mb{A}$ furnished with the multiplication
given respectively by the following formulas: 
\begin{eqnarray}
&&\hspace*{-5mm}(\varphi \bowtie u)(\psi \bowtie u{'})\nonumber\\
&&\hspace*{5mm}
=(\O ^1_{L_l}\cd \varphi \cd
\O ^5_{L_l})(\O ^2_{L_l}u_{\le0\ri _{[-1]}}\cd \psi \cd \smi (u_{\le1\ri})\O ^4_{L_l})\bowtie
\O ^3_{L_l}u_{\le0\ri _{[0]}}u{'}, \label{lgdp1}\\
&&\hspace*{-5mm}
(\v \btrl u)(\psi \btrl u{'})\nonumber\\
&&\hspace*{5mm}
=(\Omega ^1_{L_r}\cd \v \cd \Omega ^5_{L_r})
(\Omega ^2_{L_r}u_{[-1]}\cd \psi \cd \smi (u_{[0]_{\le1\ri}})\Omega ^4_{L_r})\btrl 
\Omega ^3_{L_r}u_{[0]_{\le0\ri}}u{'},\label{lgdp2}  
\end{eqnarray}
for all $\v , \psi \in {\cal A}$ and $u, u{'}\in {\mb A}$.
We write $\v \bowtie u$ (respectively $\v \btrl u$) for $\v \ot u$ considered as
an element of ${\cal A}\bowtie _{\d _l}{\mb A}$ 
(respectively ${\cal A}\btrl _{\d _r}{\mb A}$). ${\cal A}\bowtie _{\d _l}{\mb A}$ 
and ${\cal A}\btrl _{\d _r}{\mb A}$ 
are isomorphic associative algebras with unit 
$1_{\cal A}\bowtie 1_{\mb A}$ (respectively $1_{\cal A}\btrl 1_{\mb A}$),
containing ${\mb A}\cong 1_{\cal A}\bowtie {\mb A}$ 
(respectively ${\mb A}\cong 1_{\cal A}\btrl {\mb A}$)  
as unital subalgebra. These algebras are called the left 
generalized diagonal crossed products.

The right generalized diagonal 
crossed products are introduced in a similar way: denote
$$\O _{R_{l/r}}=\O ^1_{R_{l/r}}\ot \cdots \ot \O ^5_{R_{l/r}},$$
and let $\mb{A}\bowtie _{\d _l}{\cal A}$ and 
$\mb{A}\btrl _{\d _r}{\cal A}$ be the vector space $\mb{A}\ot {\cal A}$
with the following product:
\begin{eqnarray}
&&\hspace*{-5mm}
(u\bowtie \v )(u'\bowtie \psi )\nonumber\\
&&\hspace*{5mm} 
=uu'_{\le0\ri _{[0]}}\O ^3_{R_l}\bowtie (\O ^2_{R_l}
\smi (u'_{\le0\ri _{[-1]}})\cd \v \cd u'_{\le1\ri}\O ^4_{R_l})
(\O ^1_{R_l}\cd \psi \cd \O ^5_{R_l}),\label{rgdp1}\\
&&\hspace*{-5mm}
(u\btrl \v )(u'\btrl \psi )\nonumber\\
&&\hspace*{5mm}
=uu'_{[0]_{\le0\ri}}\O ^3_{R_r}\btrl (\O ^2_{R_r}\smi (u'_{[-1]})
\cd \v \cd u'_{[0]_{\le1\ri}}
\O ^4_{R_r})(\O ^1_{R_r}\cd \psi \cd \O ^5_{R_r}),\label{rgdp2}
\end{eqnarray}
for all $u, u'\in \mb{A}$ and $\v , \psi \in {\cal A}$.
$\mb{A}\bowtie _{\d _l}{\cal A}$ and 
$\mb{A}\btrl _{\d _r}{\cal A}$ are isomorphic associative algebras
with unit $1_{\mb{A}}\bowtie 1_{\cal A}$ and $1_{\mb{A}}\btrl 1_{\cal A}$,
containing $\mb{A}$ as a unital subalgebra.

As algebras, the left and right generalized crossed product algebras
are isomorphic, see \cite{hn1,bpv3}. If $H$ is a quasi-Hopf algebra then ${\cal A}=H^*$ is 
an $H$-bimodule algebra. In this particular case the left and right 
generalized diagonal crossed products are exactly the left and the right 
diagonal crossed products constructed in \cite{hn1}. In this way Hausser and Nill 
gave four explicit 
realizations of $D(H)$, the quantum double of a finite dimensional quasi-Hopf 
algebra $H$. Two of them are build on $H^*\ot H$ and the other two on $H\ot H^*$. 
All these are, as algebras, diagonal crossed products. The first two realizations 
built on $H^*\ot H$ coincide to $H^*\bowtie _{\d _l}H$ and $H^*\btrl _{\d _r}H$, 
and the last two realizations built on $H\ot H^*$ coincide to 
$H\bowtie _{\d _l}H^*$ and $H\btrl _{\d _r}H^*$.  

\begin{proposition}\prlabel{3.10}
Let $H$ be a quasi-Hopf algebra, $\mb{A}$ an $H$-bicomodule algebra 
and $C$ an $H$-bimodule coalgebra. Then $\mb{A}^1\gtl C^*\equiv 
\mb{A}\bowtie _{\d _l}C^*$ and $\mb{A}^2\gtl C^*\equiv 
\mb{A}\btrl _{\d _r}C^*$, as algebras. In particular, the algebras 
$\mb{A}^1\gtl C^*$ and $\mb{A}^2\gtl C^*$ are isomorphic to each other,  
and also to the four generalized diagonal crossed products. 
\end{proposition}

\begin{proof}
Keeping the above concepts and notations it is not hard to see that 
for an $H$-bicomodule algebra $\mb{A}$ the 
reassociators $\Phi _{\r _1}$ and $\Phi _{\r _2}$ defined in (\ref{rcaspr2}) 
and (\ref{rcaspr4}), respectively, can be rewritten as 
\begin{eqnarray*}
&&\Phi _{\r _1}=\tilde{\O }^3_{R_l}\ot (\tilde{\O }^2_{R_l}\ot 
\tilde{\O }^4_{R_l})\ot (\tilde{\O }^1_{R_l}\ot 
\tilde{\O }^5_{R_l}),\\
&&\Phi _{\r _2}=\tilde{\O }^3_{R_r}\ot (\tilde{\O }^2_{R_r}\ot 
\tilde{\O }^4_{R_r})\ot (\tilde{\O }^1_{R_r}\ot 
\tilde{\O }^5_{R_r}),
\end{eqnarray*}
where we used the notation
$$\O ^{-1}_{R_{l/r}}=\tilde{\O }^1_{R_{l/r}}\ot 
\cdots \ot \tilde{\O }^5_{R_{l/r}}.$$
 Now, if $C$ is an $H$-bimodule 
coalgebra viewed as a left $H^{\rm op}\ot H$-module coalgebra via the 
structure defined in (\ref{rmcs}) then $C^*$, the linear dual space of $C$, 
is a right $H^{\rm op}\ot H$-module algebra. The multiplication of $C^*$ is the convolution, 
that is $(c^*d^*)(c)=c^*(\una )d^*(\unb )$, the unit is $\une $ and the  
right $H^{\rm op}\ot H$-module action is given by the formula
$(c^*\cd (h\ot h'))(c)=c^*(h{'}\cd c\cd h)=(h\rh c^*\lh h')(c)$, 
for all $h, h'\in H$, $c^*, d^*\in C^*$, $c\in C$. 

Since $\mb{A}^{1/2}$ are right $H^{\rm op}\ot H$-comodule algebras and $C^*$ is a 
right $H^{\rm op}\ot H$-module coalgebra, it makes sense 
to consider the right generalized smash product algebras $\mb{A}^{1/2}\gtl C^*$
(cf. \resref{2.8} i)). It also follows easily  
from \resref{2.8} i) that the multiplication on 
$\mb{A}^1\gtl C^*$ is given by 
\begin{eqnarray*}
&&\hspace*{-1cm}
(u\gtl c^*)(u'\gtl d^*)\\
&&=uu'_{\le0\ri _{[0]}}\O ^3_{R_l}\ot 
\left((c^*\cd (\O ^2_{R_l}\smi (u'_{\le0\ri _{[-1]}})\ot u'_{\le1\ri}\O ^4_{R_l}))
(d^*\cd (\O ^1_{R_l}\ot \O ^5_{R_l}))\right)\\
&&=uu'_{\le0\ri _{[0]}}\O ^3_{R_l}\ot 
(\O ^2_{R_l}\smi (u'_{\le0\ri _{[-1]}})\rh c^*\lh u'_{\le1\ri}\O ^4_{R_l})
(\O ^1_{R_l}\rh d^*\lh \O ^5_{R_l}).
\end{eqnarray*} 

On the other hand, it is easy to see that the linear dual space $C^*$
of an $H$-bimodule coalgebra 
$C$ is an $H$-bimodule algebra. 
The multiplication of $C^*$ is the convolution, the unit is $\une $ and 
the  $H$-bimodule structure is given by the formula
$(h\rh c^*\lh h{'})(c)=c^*(h{'}\cd c\cd h)$, for all $h, h'\in H$, 
$c^*\in C^*$, $c\in C$. So we can consider the generalized right 
diagonal crossed product $\mb{A}\bowtie _{\d _l}C^*$. From (\ref{rgdp1}) 
it follows that the multiplication rule on $\mb{A}\bowtie _{\d _l}C^*$ 
coincides with the multiplication of $\mb{A}^1\gtl C^*$, hence the algebras
$\mb{A}\bowtie _{\d _l}C^*$ and $\mb{A}^1\gtl C^*$ are equal. In a similar 
way we can show the equality of the $k$-algebras
 $\mb{A}^2\gtl C^*$ and $\mb{A}\btrl _{\d _r}C^*$. 
\end{proof}

\begin{remark}
It was shown in \cite{bpv3} 
that the left generalized crossed product algebras ${\cal A}\bowtie _{\d _l}\mb{A}$ 
and ${\cal A}\btrl _{\d _r}\mb{A}$ coincide with the left generalized 
smash product algebras ${\cal A}\gsm \mb{A}_1$ and ${\cal A}\gsm \mb{A}_2$, 
respectively. The generalized smash products are made over $H\ot H^{\rm op}$ and 
by $\mb{A}_1$ and $\mb{A}_2$ we denote the left $H\ot H^{\rm op}$-comodule 
algebras structures on $\mb{A}$ defined at the end of \seref{1}. 
\end{remark}

Let $H$ be a quasi-bialgebra, $\mfA$ a right $H$-comodule algebra and 
$C$ a left $H$-module algebra. Viewing the category 
${}_{\mfA}{\cal M}(H)^C$ as a category of right-left Doi-Hopf modules,
we deduce from \thref{2.5} that ${}_{\mfA}{\cal M}(H)^C$ is isomorphic to 
the category of rational $\mfA \gtl C^*$-modules, 
$\Rat({}_{\mfA \gtl C^*}{\cal M})$, and 
%equivalent to the category 
that $\Rat({}_{\mfA \gtl C^*}{\cal M})=\sigma _{\mfA \gtl C^*}[\mfA \ot C]$.
A rational $\mfA \gtl C^*$-module 
is a left $\mfA \gtl C^*$-module $M$ such that for any $m\in M$ there 
exist two finite families $\{c_i\}_i\subseteq C$ and $\{m_i\}_i\subseteq M$ 
such that 
$$
(\mfa \gtl c^*)\cd m=c^*(c_i) (\mfa \gtl \une )\cd m_i,$$
for all $\mfa \in \mfA$ and $c^*\in C^*$.

\begin{corollary}\colabel{3.11}
Let $H$ be a quasi-Hopf algebra, $\mb{A}$ an $H$-bicomodule algebra 
and $C$ an $H$-bimodule coalgebra. Then the following assertions hold:
\begin{itemize}
\item[i)] The category of left-right Yetter-Drinfeld modules 
${}_{\mb{A}}{\cal YD}(H)^C$ is isomorphic to the category of rational 
$\mb{A}\btrl _{\d _r}C^*$-modules $\Rat({}_{\mb{A}\btrl _{\d _r}C^*}{\cal M})$, 
which is also equal to the 
category $\sigma _{\mb{A}\btrl _{\d _r}C^*}[\mb{A}\ot C]$. 
\item[ii)] If $C$ is finite dimensional then ${}_{\mb{A}}{\cal YD}(H)^C$ is 
isomorphic to the category of left $C^*\bowtie _{\d _l}\mb{A}$-modules. 
\end{itemize} 
\end{corollary}

\begin{proof}
The assertion i) follows easily from the above comments and \thref{3.8}.

ii) If $C$ is finite dimensional then 
$\Rat({}_{\mb{A}\btrl _{\d _r}C^*}{\cal M})
={}_{\mb{A}\btrl _{\d _r}C^*}{\cal M}$. Moreover, $\mb{A}\btrl _{\d _r}C^*$ 
is always isomorphic to $C^*\bowtie _{\d _l}\mb{A}$ as an algebra. We note that 
another proof of this result can be found in \cite{bpv3}. 
\end{proof}

%%%%%%%%%%%%%%%%%%%%%%%%%%%%%%%%%%%%%%%%%%%%%%%%%%%%%%%%%%%%
  
\end{document}